\documentclass[aap,seceqn,MSNbibl,citesort,dvips]{arximspdf}

\doi{10.1214/10-AAP745}
\volume{22}
\issue{1}
\pubyear{2012}
\firstpage{1}
\lastpage{28}

\makeatletter

\newproclaim{defi}{Definition}
\newproclaim{rem}{Remark}
\newtheorem{lemma}{Lemma}
\newtheorem{theo}{Theorem}
\newtheorem{corollary}{Corollary}

\newcommand{\ho}{\mathrm{hom}}
\newcommand{\R}{\mathbb{R}}
\newcommand{\Z}{\mathbb{Z}}
\newcommand{\N}{\mathbb{N}}
\newcommand{\Id}{\mathrm{Id}}
\newcommand{\med}{\int\hspace*{-11.2pt}-}
\newcommand{\calA}{\mathcal{A}_{\alpha\beta}}
\newcommand{\ee}{\mathbf{e}}

\newcommand{\eqref}[1]{(\ref{#1})}

\newcommand{\dig}[1]{\operatorname{diag} [ #1 ]}
\newcommand{\cov}[2]{\operatorname{cov} [#1;#2 ]}
\newcommand{\step}[1]{\textit{Step} #1.}

\makeatother

\begin{document}
\begin{frontmatter}

\title{An optimal error estimate in stochastic homogenization of
discrete elliptic equations\thanksref{TITL1}}
\thankstext{TITL1}{Supported in part of the Hausdorff Center for
Mathematics, Bonn, Germany.}
\runtitle{Optimal error in stochastic homogenization}

\begin{aug}
\author[A]{\fnms{Antoine} \snm{Gloria}\corref{}\ead[label=e1]{antoine.gloria@inria.fr}}
and
\author[B]{\fnms{Felix} \snm{Otto}\ead[label=e2]{otto@mis.mpg.de}}
\runauthor{A. Gloria and F. Otto}
\affiliation{INRIA Lille-Nord Europe and Max Planck Institute for \\
\quad Mathematics in the Sciences}
\address[A]{Project-team SIMPAF\\
INRIA Lille-Nord Europe\\
Parc Scientifique de la Haute Borne\\
40, avenue Halley Bat. A\\
Park Plaza\\
59650 Villeneuve d'Ascq\\
France\\
\printead{e1}} %adresu isvedimo komanda gale!
\address[B]{Max Planck Institute\\
\quad for Mathematics in the Sciences\\
Inselstr. 22\\
04103 Leipzig\\
Germany\\
\printead{e2}}
\end{aug}

% HISTORY:
\received{\smonth{2} \syear{2010}}
\revised{\smonth{8} \syear{2010}}

% ABSTRACT
%
\vspace*{24pt}
\begin{abstract}
This paper is the companion article to [\textit{Ann. Probab.} \textbf{39}
(2011) 779--856]. %\cite{Gloria-Otto-09}
We consider a discrete elliptic equation on the $d$-dimensional
lattice $\mathbb{Z}^d$ with random coefficients $A$
of the simplest type: They
are identically distributed and independent from edge to edge.
On scales large w.r.t. the lattice spacing (i.e., unity), the solution
operator is known to behave like
the solution operator of a (continuous) elliptic equation with
constant deterministic coefficients. This symmetric ``homogenized''
matrix $A_{\mathrm{hom}}=a_{\mathrm{hom}}\mathrm{Id}$ is
characterized by
$
\xi\cdot A_{\mathrm{hom}}\xi=
\langle(\xi+\nabla\phi)\cdot A(\xi+\nabla\phi)\rangle
$
for any direction $\xi\in\mathbb{R}^d$,
where the random field $\phi$ (the ``corrector'')
is the unique solution of
$
-\nabla^*\cdot A(\xi+\nabla\phi) = 0
$
in $\mathbb{Z}^d$ such that $\phi(0)=0$, $\nabla\phi$ is stationary
and $\langle\nabla\phi\rangle=0$,
$\langle\cdot\rangle$ denoting the ensemble average (or expectation).

In order to approximate the homogenized coefficients $A_{\mathrm
{hom}}$, the corrector
problem is usually solved in a box $Q_L=[-L,L)^d$ of size $2L$ with
periodic boundary
conditions, and the space averaged energy on $Q_L$ defines an
approximation $A_L$ of
$A_{\mathrm{hom}}$. Although the statistics is modified (independence
is replaced
by periodic correlations) and the ensemble average is replaced by a
space average,
the approximation $A_L$ converges almost surely to $A_{\mathrm{hom}}$
as $L \uparrow\infty$.
In this paper, we give estimates on both errors. To be more precise, we
do not consider
periodic boundary conditions on a box of size $2L$, but replace the
elliptic operator by
$
T^{-1}-\nabla^*\cdot A\nabla
$
with (typically) $T\sim L^2$, as standard in the homogenization literature.
We then replace the ensemble average by a space average on $Q_L$,
and estimate the overall error on the homogenized coefficients in terms
of $L$ and $T$.
% \medskip
% To improve the convergence, one may also consider $N$ independent
% realizations of $A_L$ and take their arithmetic mean.
% A natural optimization problem consists in properly choosing $L$ and
%$N$ in order to
% reduce the error at given computational complexity.
% The analysis we develop in this article, combined with our variance
%estimate in
% \cite{Gloria-Otto-09}, allows us to give a first piece of answer to
%this question.
\end{abstract}

% KEYWORDS
%

\end{frontmatter}

%s1 ###
\section{Introduction}\label{sec:intro}

%s1.1 ###
\subsection{Motivation}

In this article, we continue the analysis we began in \cite{Gloria-Otto-09}
on stochastic homogenization of discrete elliptic equations. More precisely,
we consider real functions $u$ of the sites $x$
in a $d$-dimensional Cartesian lattice~$\mathbb{Z}^d$.
Every edge $e$ of the lattice is endowed with a ``conductivity''
$a(e)>0$. This defines a discrete elliptic differential operator
$-\nabla^*\cdot A\nabla$ via
\[
-\nabla^*\cdot(A\nabla u)(x) := \sum_{y\in\Z^d,|x-y|=1} a(e)\bigl
(u(x)-u(y)\bigr),
\]
where the sum is over the $2d$ sites $y$ which are connected by
an edge $e=[x,y]$ to the site $x$. It is sometimes more convenient
to think in terms of the associated Dirichlet form, that is,\
\begin{eqnarray*}
\sum\nabla v\cdot A\nabla u
&:=&
\sum_{x\in\mathbb{Z}^d}v(x) \bigl(-\nabla^*\cdot(A\nabla
u)(x)\bigr)\\
&=&
\sum_{e}\bigl(v(x)-v(y)\bigr)a(e)\bigl(u(x)-u(y)\bigr),
\end{eqnarray*}
where the last sum is over all edges $e$ and $(x,y)$ denotes the two sites
connected by $e$, that is, $e=[x,y]=[y,x]$ (with the convention that an
edge is not oriented).
We assume the conductivities
$a$ to be uniformly elliptic in the sense of
\[
\alpha\le a(e) \le\beta\qquad\mbox{for all edges } e
\]
for some fixed constants $0<\alpha\le\beta<\infty$.

We are interested in random coefficients. To fix ideas,
we consider the simplest situation possible:
\[
\{a(e)\}_{e} \mbox{ are independently and identically distributed
(i.i.d.)}.
\]
Hence, the statistics are described by a distribution on the finite
interval $[\alpha,\beta]$. We'd like to see this discrete elliptic operator
with random coefficients as a good model problem for
continuum elliptic operators with random coefficients of correlation
length unity.

Classical results in stochastic homogenization of linear
elliptic equations
(see \cite{Kozlov-79} and \cite{Papanicolaou-Varadhan-79}
for the continuous case, and \cite{Kunnemann-83} and \cite{Kozlov-87}
for the discrete case)
state that there exist \textit{homogeneous and deterministic} coefficients
$A_{\ho}$
such that the solution operator of the continuum differential operator
$-\nabla\cdot A_{\ho}\nabla$
describes the large scale behavior of the solution operator
of the discrete differential operator $-\nabla^*\cdot A\nabla$.
As a by product of this homogenization result,
one obtains a characterization of the homogenized
coefficients $A_{\ho}$: It is shown that for
every direction $\xi\in\mathbb{R}^d$, there exists a unique
scalar field~$\phi$ such that $\nabla\phi$ is stationary
[stationarity means that the fields
$\nabla\phi(\cdot)$ and $\nabla\phi(\cdot+z)$ have the same\vadjust{\goodbreak} statistics
for all shifts $z\in\mathbb{Z}^d$] and $\langle\nabla\phi \rangle=0$,
solving the equation
%
%e1.1 ###
\begin{equation}\label{PV2}
-\nabla^*\cdot\bigl(A(\xi+\nabla\phi)\bigr) = 0 \qquad\mbox{in
} \mathbb{Z}^d
\end{equation}
and normalized by $\phi(0)=0$. As in periodic homogenization,
the function $\mathbb{Z}^d\ni x\mapsto\xi\cdot x+\phi(x)$
can be seen as the $A$-harmonic function which macroscopically
behaves as the affine function $\mathbb{Z}^d\ni x\mapsto\xi\cdot x$.
With this
``corrector''~$\phi$, the homogenized coefficients $A_{\ho}$
(which in general form a symmetric matrix and for our simple statistics
in fact a multiple of the identity: $A_{\ho}=a_{\ho}\Id$)
can be characterized as follows:
%
%e1.2 ###
\begin{equation}\label{PV1}
%a_{\ho}\xi&=&\langle(a(\xi+\nabla\phi)\rangle, \mbox{or
%equivalently}\\
\xi\cdot A_{\ho}\xi=
\langle(\xi+\nabla\phi)\cdot A(\xi+\nabla\phi)\rangle.
\end{equation}
Since
%the vector field $a(\xi+\nabla\phi)$ and
the scalar field
$(\xi+\nabla\phi)\cdot A(\xi+\nabla\phi)$ is stationary, it does
not matter (in terms of the distribution) at which site $x$ it
is evaluated in the formula~(\ref{PV1}), so that we suppress
the argument $x$ in our notation.

When one is interested in explicit values for $A_\ho$, one
has to solve \eqref{PV2}.
Since this is not possible in practice, one has to make approximations.
For a discussion of the literature on error estimates, in particular
the pertinent work by Yurinskii
\cite{Yurinskii-86} and Naddaf and Spencer \cite{Naddaf-Spencer-98},
we refer to \cite{Gloria-Otto-09}, Section~1.2.
A standard approach used in practice consists in solving \eqref{PV2}
in a~box $Q_L=[-L,L)^d\cap\Z^d$ with periodic boundary conditions
%
%e1.3 ###
\begin{equation}\label{PV2-b}
-\nabla^*\cdot\bigl(A(\xi+\nabla\phi_{L,\#})\bigr) = 0\qquad
\mbox{in } Q_L,
\end{equation}
and replacing \eqref{PV1} by a space average
%
%e1.4 ###
\begin{equation}\label{PV1-b}
\xi\cdot A_{L,\#}\xi=
\med_{Q_L}(\xi+\nabla\phi_{L,\#})\cdot A(\xi+\nabla\phi
_{L,\#})\,dx.
\end{equation}
Such an approach is consistent in the sense that
\[
\lim_{L\to\infty} A_{L,\#} = A_\ho
\]
almost surely, as proved, for instance, in \cite{Bourgeat-04} for the
continuous case, and in~\cite{Caputo-Ioffe-03} for the discrete case.
Numerical experiments tend to show that the use of periodic boundary
conditions gives better results than other choices
such as homogeneous Dirichlet boundary conditions, see \cite{E-Yue-07}.

An important question for practical purposes is to quantify
the dependence of the
error $\langle|A_\ho-A_{L,\#}|^2 \rangle^{1/2}$ in terms of $L$.
Let us give another interpretation of \eqref{PV2-b}: This equation on
$Q_L$ is equivalent to \eqref{PV2}
on $\Z^d$ with a modified conductivity matrix $\tilde A_L$, that is
the periodization of $A_{|Q_L}$ on $\Z^d$.
Doing this, we have replaced independent coefficients $A$ by
$Q_L$-periodically correlated coefficients $\tilde A$.
Since $A$ and $\tilde A$ are not jointly stationary (see
Definition~\ref{def:stationary}), it may be
difficult to compare $\nabla\phi$ to $\nabla\phi_{L,\#}$.
To circumvent this difficulty, and following the route of\vadjust{\goodbreak} \cite
{Papanicolaou-Varadhan-79,Kunnemann-83,Yurinskii-86} and \cite
{Naddaf-Spencer-98}, and as in \cite
{Gloria-Otto-09},
we slightly depart from \eqref{PV2-b} by introducing a zero-order term
in \eqref{PV2}:
%
%e1.5 ###
\begin{equation}\label{PV2-c}
T^{-1}\phi_T-\nabla^*\cdot\bigl(A(\xi+\nabla\phi_T)\bigr) = 0
\qquad\mbox
{in } \mathbb{Z}^d.
\end{equation}
As for the periodization, this localizes the dependence of $\phi_T(z)$
upon $A(z')$ to those points $z'\in\Z^d$ such
that $|z-z'| \lesssim\sqrt{T}$ (at first order). Yet, unlike the
periodization, $\nabla\phi_T$ and $\nabla\phi$ are
jointly stationary.
In terms of random walk interpretation, the lifetime of the random
walker is of order $T$,
and the distance to the origin of order $\sqrt{T}$.
Hence, up to taking $T \sim L^2$, in first approximation, the function
${\phi_T}_{|Q_L}$ only depends on the
coefficients $A(z)$ for $z\in Q_L$, as it is the case for $\phi_{L,\#}$.

We'd like to view ${\phi_T}_{|Q_L}$ as a variant of $\phi
_{L,\#}$ which is convenient for our analysis.
We then define
%
%e1.6 ###
\begin{equation}\label{PV1-c}
\xi\cdot A_{T,L}\xi=
\int_{\Z^d}(\xi+\nabla\phi_{T})\cdot A(\xi+\nabla\phi
_{T})\eta_L\,dx,
\end{equation}
where $\eta_L$ is a smooth mask with unit mass and support $Q_L$.
The aim of this paper is to determine the scaling of the error
$\langle|A_\ho-A_{T,L}|^2 \rangle^{1/2}$ in terms of $L$ and $T$.
Eventually this will allow us to make a reasonable choice for $T$ and
$L$ at fixed computational complexity.

%s1.2 ###
\subsection{Informal statement of the results}

$\!\!\!$When approximating $A_\ho$ by~$A_{T,L}$, we make two types of errors:
A ``systematic error'' and a ``random error.''
In particular, as shown in \cite{Gloria-Otto-09},
\[
\langle(\xi\cdot A_\ho\xi-\xi\cdot A_{T,L}\xi)^2 \rangle = \bigl
(\xi
\cdot(A_\ho-A_T)\xi\bigr)^2+\operatorname{var} [A_{T,L} ].
\]
The first term is the square of the systematic error (see \cite
{Gloria-Otto-09}, (1.10))
\begin{eqnarray}\label{eq:err-syst}
\mathrm{Error}_{\mathrm{sys}}(T)&:=&|\langle(\xi+\nabla\phi
_T)\cdot A (\xi+\nabla\phi_T) \rangle-\langle(\xi+\nabla\phi
)\cdot A (\xi+\nabla\phi) \rangle| \nonumber\hspace*{-35pt}\\[-8pt]\\[-8pt]
&\hspace*{3pt}=&\langle(\nabla\phi_T-\nabla\phi)\cdot A (\nabla
\phi_T-\nabla \phi) \rangle.\nonumber\hspace*{-35pt}
\end{eqnarray}
It measures the fact that the coefficient $a(e)$ at bond $e$ does (up
to exponentially small terms) not influence $\phi_T(x)$ if
$|x-e|\gg\sqrt{T}$. This error vanishes for $T=L^2\uparrow\infty$.
The second term is the square of the random error,
%
%e1.7 ###
\begin{eqnarray}\label{eq:err-rand}
\mathrm{Error}_{\mathrm{rand}}(T,L)
&=&\operatorname{var}\biggl[\int_{\Z^d} (\xi+\nabla\phi
_{T})\cdot A (\xi+\nabla\phi
_{T})\eta_L\,dx \biggr]^{1/2}.
\end{eqnarray}
It measures the fluctuations of the energy density. This error vanishes
as $L\uparrow\infty$.

In \cite{Gloria-Otto-09}, Theorem~1, we have proved that
%
%e1.8 ###
\begin{equation}\label{eq:intro-Thm-var}
\qquad\operatorname{var}\biggl[\int_{\Z^d} (\xi+\nabla\phi
_{T})\cdot A
(\xi+\nabla\phi
_{T})\eta_L \,dx\biggr]^{1/2}
\lesssim
\cases{
d=2,&\quad$L^{-1}\ln^qT$, \cr
d>2,&\quad$L^{-d/2}$,
}\hspace*{-30pt}
\end{equation}
for some $q$ depending only on $\alpha,\beta$, where ``$\lesssim$''
stands for ``$\leq$'' up to a~multiplicative
constant depending only on $\alpha,\beta$ and $d$.
We have also identified the\vadjust{\goodbreak} systematic error in the limit of vanishing
conductivity contrast, that is,
$1-\beta/\alpha\ll1$, and found
\[
\mathrm{Error}_{\mathrm{sys}}(T) \sim
\cases{
d=2,&\quad$T^{-1}$, \cr
d=3,&\quad$T^{-3/2}$,\cr
d=4,&\quad$T^{-2}\ln T$,\cr
d>4,&\quad$T^{-2}$,
}
\]
where ``$\sim$'' means that both terms have the same scaling (in $T$).
In this paper, we shall actually prove that for general $\alpha$ and
$\beta$ (see Theorem \ref{th1})
%
%e1.9 ###
\begin{equation}\label{eq:error-sys}
\mathrm{Error}_{\mathrm{sys}}(T) \lesssim
\cases{
d=2,&\quad$T^{-1}\ln^qT$, \cr
d=3,&\quad$T^{-3/2}$,\cr
d=4,&\quad$T^{-2}\ln T$,\cr
d>4,&\quad$T^{-2}$,
}
\end{equation}
where there is a logarithmic correction for $d=2$ when compared to the
vanishing conductivity asymptotics.

Assuming that $\phi_T$ can be well approximated on domains
of size $L$ if we choose $T\sim L^{2}$,
the combination of \eqref{eq:error-sys} and \eqref{eq:intro-Thm-var} yields
\[
\langle|A_\ho-A_{T,L}|^2 \rangle^{1/2} \lesssim
\cases{
d=2,&\quad$L^{-1}\ln^qL $,\cr
2<d \leq7,&\quad$L^{-d/2}$, \cr
d=8,&\quad$L^{-4}\ln L $,\cr
d>8,&\quad$L^{-4}$.
}
\]
Hence, the numerical strategy converges at the rate of the central
limit theorem for $2\leq d \leq8$ (up to logarithmic corrections for
$d=2$ and $d=8$).

Up to dimension 4, the systematic error for $T\sim L^2$
scales as the square of the random error.
In particular, this leaves room for the choice $T$.
If we take $T\sim L$, then the systematic error is of the same order as
the random error.
What we have gained is that $\phi_T$ can now be well-approximated on
domains of size $R\sim\sqrt{T} \sim\sqrt{L}$,
and not only $L$.
Note also that the random error is unchanged if instead of taking the
average of one realization of $\phi_T$
on~$Q_L$ (with the mask $\mu_L$) we take the empirical average of the
averages of~$N$ independent realizations
of $\phi_T$ on a domain $Q_{L/N^{1/d}}$ (with the according mask $\mu
_{L/N^{1/d}}$).
Hence, since $\phi_T$ can be well-approximated on domains of size
$R\gtrsim\sqrt{L}$,
considering $N=1$ realization of $\phi_T$ approximated on $Q_L$ or
$N=\sqrt{L^d}$ independent realizations of $\phi_T$
approximated on $Q_{\sqrt{L}}$ yields
the same scaling for the error between the homogenized coefficients and
their approximations.
Since the computational cost of solving a linear problem is superlinear
in the number of unknowns, it seems best
to choose $N$ as large as possible, and therefore taking $N=\sqrt
{L^d}$ seems a reasonable strategy at first order.
Yet, we do not make precise in this paper the relation between~$R$ and
$\sqrt{L}$ in terms of absolute values (we only consider the scaling),
which may make the optimal choice for $N$ more subtle in practice than
this general principle.
A complete numerical analysis of the numerical method (including the
influence of $R$ and the optimization of $N$) will be presented in
\cite{Gloria-10}.\vadjust{\goodbreak}

We conclude this introduction by mentioning the very recent
contribution~\cite{Mourrat-10} by Mourrat.
The equation under investigation is the same as above, namely a
discrete elliptic equation on $\Z^d$ with i.i.d.
coefficients. The object under study is the spectral measure associated
with the generator of the environment
viewed by the particle. Without entering into details, there exists
some nonnegative measure $e_{\mathfrak{d}}$ associated with
the elliptic operator and direction $\xi\in\R^d$, such that the
homogenized coefficient is given by
\[
\xi\cdot A_\ho\xi= \langle\xi\cdot A\xi \rangle-\int_{\R
^+}\frac
{1}{\lambda}\,d e_{\mathfrak{d}}(\lambda).
\]
As recalled in \cite{Mourrat-10}, we also have
%
%e1.10 ###
\begin{equation}\label{intro:spectral}
\xi\cdot A_T \xi= \langle\xi\cdot A\xi \rangle-\int_{\R^+}\frac
{\lambda-2/T}{(1/T+\lambda)^2}\,d e_{\mathfrak{d}}(\lambda).
\end{equation}
In particular, the systematic error can be written as
\[
\xi\cdot A_T \xi-\xi\cdot A_\ho\xi= \frac{1}{T^2}\int_{\R
^+}\frac{1}{\lambda(1/T+\lambda)^2}\,d e_{\mathfrak{d}}(\lambda),
\]
so that information on the scaling of the systematic error in terms of
$T$ yields information on
the spectral behavior and conversely.
The interplay between the strategy used in the present paper and the
spectral measure is further investigated
by Mourrat and the first author in \cite{Gloria-Mourrat-10}.
In what follows, we do not make use of the spectral measure, which
makes our approach self-contained.

The article is organized as follows: In Section~\ref
{sec:main-res}, we introduce the general
framework and state the main results of this paper, that is, the
systematic error actually scales as in \eqref{eq:error-sys}.
The last two sections are dedicated to its proof.

Throughout the paper, we make use of the following notation:
\begin{itemize}
\item$d\geq2$ is the dimension;
\item$\int_{\Z^d} \,dx$ denotes the sum over $x\in\Z^d$, and $\int
_{D}\,dx$ denotes the sum over $x\in\Z^d$ such that $x\in D$, $D$
subset of $\R^d$;
\item$\langle\cdot \rangle$ is the ensemble average, or
equivalently the
expectation in the underlying probability space;
\item$\operatorname{var} [\cdot]$ is the variance associated with
the ensemble average;
\item$\cov{\cdot}{\cdot}$ is the covariance associated with the
ensemble average;
\item$\lesssim$ and $\gtrsim$ stand for $\leq$ and $\geq$ up to a
multiplicative constant which only depends on the dimension $d$ and the
constants $\alpha,\beta$ (see Definition~\ref{def:alpha-beta} below)
if not otherwise stated;
\item when both $\lesssim$ and $\gtrsim$ hold, we simply write $\sim$;
\item we use $\gg$ instead of $\gtrsim$ when the multiplicative
constant is (much) larger than~$1$;
\item$(\ee_1,\dots,\ee_d)$ denotes the canonical basis of $\Z^d$.
\end{itemize}
%

%%%%%%%%%%%%%%%%%%%%%%%%%%%%%%%%%%%%%%%%%%%%%%%%%%%%%%%%%%
%%%%%%%%%%%%%%%%%%%%%%%%%%%%%%%%%%%%%%%%%%%%%%%%%%%%%%%%%%

%s2 ###
\section{Main result}\label{sec:main-res}

%s2.1 ###
\subsection{General framework}

\begin{defi}\label{def:alpha-beta}
We say that $a$ is a conductivity function if there exist $0<\alpha
\leq\beta<\infty$ such that
for every edge $e$ of $\Z^d$, one has $a(e)\in[\alpha,\beta]$.
We denote by $\calA$ the set of such conductivity functions.
\end{defi}

\begin{defi}
The elliptic operator $L\dvtx L^2_{\mathrm{loc}}(\Z^d)\to
L^2_{\mathrm
{loc}}(\Z^d)$, $u\mapsto Lu$ associated with a conductivity function
$a\in\calA$ is defined for all~$x\in\Z^d$~by
%
%e2.1 ###
\begin{equation}\label{eq:def-elliptic}
(Lu)(x)=-\nabla^*\cdot A(x)\nabla u(x),
\end{equation}
where
\[
\nabla u(x):= \left[
\matrix{
u(x+\ee_1)-u(x) \cr
\vdots\cr
u(x+\ee_d)-u(x)
}
\right],
\qquad
\nabla^* u(x):= \left[
\matrix{
u(x)-u(x-\ee_1) \cr
\vdots\cr
u(x)-u(x-\ee_d)
}
\right]
\]
and
\[
A(x):=\dig{a(e_1),\dots,a(e_d)},
\]
$e_1=[x,x+\ee_1],\dots,e_d=[x,x+\ee_d]$.
\end{defi}

We now turn to the definition of the statistics of the
conductivity function.

\begin{defi}
$\!\!\!$A conductivity function is said to be independent and~iden\-tically
distributed (i.i.d.) if the coefficients $a(e)$ are i.i.d. random variables.
\end{defi}
\begin{defi}\label{def:stationary}
The conductivity matrix $A$ is obviously stationary in the sense that
for all $z\in\Z^d$,
$A(\cdot+z)$ and $A(\cdot)$ have the same statistics, so that for all
$x,z\in\Z^d$,
\[
\langle A(x+z) \rangle=\langle A(x) \rangle.
\]
Therefore, any translation invariant function of $A$, such as the modified
corrector $\phi_T$ (see Lemma~\ref{lem:app-corr}), is jointly
stationary with $A$.
In particular, not only are $\phi_T$ and its gradient $\nabla\phi_T$
stationary, but also
any function of $A$,~$\phi_T$ and~$\nabla\phi_T$. A useful such
example is the energy density
\mbox{$(\xi\,{+}\,\nabla\phi_T)\,{\cdot}\,A (\xi\,{+}\,\nabla\phi_T)$}, which is
stationary by joint stationarity of $A$ and $\nabla\phi_T$.

Another translation invariant function of $A$ is the Green
functions $G_T$ of Definition~\ref{def:Green}.
In this case, stationarity means that $G_T(\cdot+z,\cdot+z)$ has the
same statistics as $G_T(\cdot,\cdot)$ for all
$z\in\Z^d$, so that in particular, for all~$x,y,z\in\Z^d$,
\[
\langle G_T(x+z,y+z) \rangle=\langle G_T(x,y) \rangle.
\]
\end{defi}
\begin{lemma}[(Corrector; \cite{Kunnemann-83}, Theorem~3)] \label{lem:corr}
Let $a\in\calA$ be an i.i.d. conductivity function, then for all $\xi
\in\R^d$, there exists a unique random function~$\phi\dvtx\Z^d\to
\R$
which satisfies the corrector equation
%
%e2.2 ###
\begin{equation}\label{eq:corr}
-\nabla^*\cdot A(x) \bigl( \xi+\nabla\phi(x) \bigr)=0\qquad
\mbox{in }\Z^d,\vadjust{\goodbreak}
\end{equation}
and such that $\phi(0)=0$, $\nabla\phi$ is stationary and $\langle
\nabla\phi \rangle=0$.
In addition,\break \mbox{$\langle|\nabla\phi|^2 \rangle\lesssim|\xi|^2$}.
\end{lemma}

We also define an ``approximation'' of the corrector as follows.
\begin{lemma}[(Approximate corrector; \cite{Kunnemann-83}, Proof of
Theorem~3)]\label{lem:app-corr}
Let \mbox{$a\in\calA$} be an i.i.d. conductivity function, then for all
$T>0$ and $\xi\in\R^d$, there exists a unique stationary random
function $\phi_{T}\dvtx\Z^d\to\R$ which satisfies the ``approximate''
corrector equation
%
%e2.3 ###
\begin{equation}\label{eq:app-corr}
T^{-1}\phi_{T}(x)-\nabla^*\cdot A(x) \bigl( \xi+\nabla\phi_{T}(x)
\bigr)=0\qquad
\mbox{in }\Z^d,
\end{equation}
and such that $\langle\phi_{T} \rangle=0$.
In addition, $T^{-1} \langle\phi_{T}^2 \rangle+\langle{|\nabla
\phi _{T}|}^2 \rangle\lesssim|\xi|^2$.
\end{lemma}
\begin{defi}[(Homogenized coefficients)]
Let $a\in\calA$ be an i.i.d. conductivity function and let $\xi\in
\R^d$ and $\phi$ be as in Lemma~\ref{lem:corr}.
We define the homogenized $d\times d$-matrix $A_\ho$ as
%
%e2.4 ###
\begin{equation}\label{eq:fo-homog}
\xi\cdot A_\ho\xi= \langle{(\xi+\nabla\phi)\cdot A(\xi +\nabla
\phi)(0)} \rangle.
\end{equation}
\end{defi}

Note that (\ref{eq:fo-homog}) fully characterizes $A_\ho$
since $A_\ho$ is a symmetric matrix (it is actually of the form $a_\ho
\Id$ for an i.i.d. conductivity function).

%s2.2 ###
\subsection{Statement of the main results}

The main result of the article is the following estimate of the
systematic error introduced in Section~\ref{sec:intro}.
\begin{theo}\label{th1}
Let $a\in\calA$ be an i.i.d. conductivity function, and let $\phi_T$
denote the approximate corrector associated with the conductivity
function~$a$ and direction $\xi\in\R^d$, $|\xi|= 1$. We then define
for all $T \gg1$ the symmetric matrix $A_{T}$ characterized by
%
%e2.5 ###
\begin{eqnarray}\label{eq:def-AT}
\xi\cdot A_{T}\xi:=\langle(\xi+\nabla\phi_T)\cdot A(\xi+\nabla
\phi_T) \rangle.
\end{eqnarray}
Then, there exists an exponent $q>0$ depending only on $\alpha,\beta$
such that
\begin{eqnarray} \label{eq:estim-err-hom}
d&=&2\mbox{:}\qquad |A_\ho-A_T| \lesssim T^{-1}(\ln T)^q, \nonumber\\
d&=&3\mbox{:}\qquad |A_\ho-A_T| \lesssim T^{-3/2},\nonumber\\[-8pt]\\[-8pt]
d&=&4\mbox{:}\qquad |A_\ho-A_T| \lesssim T^{-2}\ln T,\nonumber\\
d&>&4\mbox{:}\qquad |A_\ho-A_T| \lesssim T^{-2}.\nonumber
\end{eqnarray}
\end{theo}

As a by-product of the proof of Theorem~\ref{th1}, we obtain
the following corollary.

\begin{corollary}\label{prop1}
Let $a\in\calA$ be an i.i.d. conductivity function, $d>2$, $T>0$, and
let $\phi_T$ and $\tilde\phi$
denote the approximate corrector and stationary corrector\vadjust{\goodbreak}
(see \cite{Gloria-Otto-09}, Corollary~1) associated with the
conductivity function $a$ and direction $\xi\in\R^d$, $|\xi|= 1$,
respectively.
Then
%
%e2.6 ###
\begin{equation}\label{eq:prop1}
T^{-1}\langle(\phi_T-\tilde\phi)^2 \rangle +\langle|\nabla\phi
_T-\nabla \tilde\phi|^2 \rangle \lesssim\cases{
d=3,&\quad$T^{-3/2}$, \vspace*{-1pt}\cr
d=4,&\quad$T^{-2}\ln T$, \vspace*{-1pt}\cr
d>4,&\quad$T^{-2}$.
}
\end{equation}
In particular,
\[
\lim_{T\to\infty} \bigl(\langle(\phi_T-\tilde\phi)^2 \rangle
+\langle|\nabla \phi_T-\nabla\tilde\phi|^2 \rangle \bigr) = 0.\vspace*{-2pt}
\]
\end{corollary}

This corollary gives a full characterization of the
convergence of the regularized corrector to the exact corrector
for $d>2$.\vspace*{-2pt}

\begin{rem}
Note that the definition \eqref{eq:def-AT} of $A_T$ does not include
the zero-order term $T^{-1}\langle\phi_T^2 \rangle$, so that $\xi
\cdot
A_T\xi$ does not
coincide with the energy associated with the equation. Surprisingly,
the addition of the zero-order term
in the definition of $A_T$ would make the estimate \eqref
{eq:estim-err-hom} saturate at $T^{-1}$
for $d>2$.\vspace*{-2pt}
\end{rem}

\begin{rem}
For $d=2$, although we lose control of $\phi_T$ we may still quantify
the rate of convergence of $\nabla\phi_T$ to $\nabla\phi$,
the gradient of the corrector of Definition~\ref{lem:corr}. In
particular, \eqref{eq:prop1} is replaced by
\[
\langle|\nabla\phi_T-\nabla\phi|^2 \rangle \lesssim T^{-1}\ln^qT
\]
for some $q>0$ depending only on $\alpha,\beta$.\vspace*{-2pt}
\end{rem}
%

%s2.3 ###
\subsection{Auxiliary lemmas}

In order to prove Theorem~\ref{th1} and Corollary~\ref{prop1}, we
need three auxiliary lemmas in addition to the results of \cite
{Gloria-Otto-09}: The first
one is a~covariance estimate very similar to the variance estimate in
\cite{Gloria-Otto-09}, Lemma~2.3,
the next one is a refined version of the decay estimates of \cite
{Gloria-Otto-09}, Lemma~2.8,
whereas the last one is a generalization of the convolution estimate
of \cite{Gloria-Otto-09}, Lemma~2.10.\vspace*{-2pt}

\begin{lemma}[(Covariance estimate)]\label{lem:cov}
Let $a=\{a_i\}_{i\in\mathbb{N}}$ be a sequence of i.i.d. random
variables with range $[\alpha,\beta]$.
Let $X$ and $Y$ be two Borel measurable functions of $a\in\mathbb
{R}^\mathbb{N}$ (i.e., measurable
w.r.t. the smallest $\sigma$-algebra on $\mathbb{R}^\mathbb{N}$
for which all coordinate functions $\mathbb{R}^\mathbb{N}\ni a\mapsto
a_i\in\mathbb{R}$
are Borel measurable, cf.\ \cite{Klenke-08}, Definition 14.4).

Then we have
%
%e2.7 ###
\begin{equation}\label{eq:covar-estim}
\cov{X}{Y} \le
\sum_{i=1}^\infty\biggl\langle\sup_{a_i} \biggl|\frac{\partial
X}{\partial a_i} \biggr|^2 \biggr\rangle^{1/2}
\biggl\langle\sup_{a_i} \biggl|\frac{\partial Y}{\partial a_i} \biggr|^2 \biggr\rangle^{1/2}
\operatorname{var} [a_1 ],
\end{equation}
where $\sup_{a_i} |\frac{\partial Z}{\partial a_i} |$ denotes
the supremum of the modulus of
the $i$th partial derivative
\[
\frac{\partial Z}{\partial a_i}(a_1,\ldots
,a_{i-1},a_i,a_{i+1},\ldots)
\]
of $Z$
with respect to the variable\vadjust{\goodbreak} $a_i\in[\alpha,\beta]$, for $Z=X,Y$.
\end{lemma}

The proof of this lemma is standard.
As for \cite{Gloria-Otto-09}, Lemma~2.3, it relies on a~martingale
difference decomposition.

We define discrete Green's functions in the following definition.
\begin{defi}[(Discrete Green's function)]\label{def:Green}
Let $d\geq2$.
For all $T>0$, the Green function $G_T\dvtx\calA\times\Z^d\times\Z
^d\to\Z^d,(a,x,y)\mapsto G_T(x,y;a)$
associated with the conductivity function $a$ is defined for all $y\in
\Z^d$ and $a\in\calA$ as the unique solution $G(\cdot,y;a)\in
L^2(\Z^d)$ to
%
%e2.8 ###
\begin{eqnarray} \label{eq:disc-Green}
\qquad\quad \int_{\Z^d}T^{-1}G_T(x,y;a)v(x)\,dx+\int_{\Z^d}\nabla v(x)\cdot
A(x)\nabla_x G_T(x,y;a)\,dx=v(y)\nonumber\\[-8pt]\\[-8pt]
\eqntext{\forall v\in L^2(\Z^d),}
\end{eqnarray}
where $A$ is as in~(\ref{eq:def-elliptic}).
\end{defi}

Throughout this paper, when no confusion occurs, we use the
shorthand notation $G_T(x,y)$ for $G_T(x,y;a)$.
We need a decay of the Green function~$G_T(x,y)$ and its (discrete)
gradient $\nabla_x G_T(x,y)$ in $|x-y|\gg1$
that is \textit{uniform} in $a$ but nevertheless coincides (in terms of scaling)
with the decay of the \textit{constant-coefficient} Green function.
The constant-coefficient
Green function in the continuous case is known to decay as
\begin{eqnarray*}
&&|x-y|^{2-d}\exp\biggl(-\mbox{const.} \frac{|x-y|}{\sqrt{T}}\biggr)\qquad\hspace*{10pt} \mbox
{for }d>2 \quad \mbox{and}\\
&&\biggl(\ln\frac{\sqrt{T}}{|x-y|}\biggr)\exp\biggl(-\mbox{const.} \frac{|x-y|}{\sqrt
{T}}\biggr)\qquad \mbox{for }d=2;
\end{eqnarray*}
its gradient decays as the first derivative of these expressions. Note
the cross-over
of the decay at distances $|x-y|$ of the order of the intrinsic length scale
$\sqrt{T}\gg1$ from algebraic (or logarithmic in case of $d=2$) to
exponential.\vadjust{\goodbreak}

In the class of $a$-uniform estimates, these decay properties
survive as
\textit{pointwise} in $(x,y)$ estimates on the level of the discrete Green
function~$G_T(x,y)$
itself, but only as \textit{averaged} estimates on the level of its
discrete gradient~$\nabla_x G_T(x,y)$.
More precisely, $\nabla_x G_T(x,y)$ has to be averaged in $x$ on
dyadic annuli centered at $x=y$. It will be important that
the average can be (at least slightly) stronger than a \textit{square} average
(see \cite{Gloria-Otto-09}, Lemma~2.9).
On the other hand, we do not need the exponential decay: Super
algebraic decay
is sufficient for our purposes.
\begin{lemma}[(Pointwise decay estimate on $G_T$)]\label{lem:ptwise-Green-decay}
Let $a\in\calA$, and $G_T$ be the associated Green function.
For $d>2$, we have for all $k>0$, and all~$x,y\in\Z^d$
%
%e2.9 ###
\begin{equation}\label{eq:def_gT_d>2}
G_T(x,y) \lesssim(1+|x-y|)^{2-d}\min\biggl\{1,\biggl(\frac{|x-y|}{\sqrt
{T}}\biggr)^{-k}\biggr\},
\end{equation}
where the constant in ``$\lesssim$'' depends on $k$.
For $d=2$, we have for all $k>0$
%
%e2.10 ###
\begin{equation}\label{eq:def_gT_d=2}
G_T(x,y) \lesssim\left\{
\matrix{
\displaystyle \ln\biggl(\frac{\sqrt{T}}{1+|x-y|}\biggr)\mbox{ for }|x-y|\ll\sqrt{T}\vspace*{2pt}\cr
\displaystyle \biggl(\frac{|x-y|}{\sqrt{T}}\biggr)^{-k}\mbox{ for }|x-y|\gtrsim\sqrt{T}
}
\right\},
\end{equation}
where the constant in ``$\lesssim$'' depends on $k$.
\end{lemma}

Finally, for the proof of Theorem~\ref{th1}, we need to know
that also the \textit{convolution} of the gradient of the Green's
function with itself decays at the optimal rate, that is, with the
following lemma.
\begin{lemma}[(Convolution estimate)]\label{b-lem:dyad}
Let $h_T,g_T\dvtx\Z^d\to\R^+$ satisfy the following properties.

\emph{Assumptions on $h_T$ [estimate of $|\nabla
_xG_T(y+z,y)|$]:} For all $R\gg1$ and $T>0$,
%
%e2.12 ###
%e2.11 ###
\begin{eqnarray}
d&>&2\mbox{:}\qquad  \int_{R< |z|\leq2R}h_T(z)^2\,dz \lesssim R^{2-d}, \label
{b-eq:assump-h1-d>2}\\
d&=&2\mbox{:}\qquad \int_{R< |z| \leq2R}h_T(z)^2\,dz\lesssim\min\bigl\{1,\sqrt
{T}R^{-1}\bigr\}^2, \label{b-eq:assump-h1-d=2}
\end{eqnarray}
and for $R\sim1$
%
%e2.13 ###
\begin{equation}\label{b-eq:assump-h4-d>=2}
d\geq2\mbox{:} \qquad  \int_{ |z|\leq R}h_T(z)^2\,dz \lesssim1 .
\end{equation}

\emph{Assumptions on $g_T$ [estimate of $G_T(y+z,y)$]:} For
$d>2$, and for all \mbox{$z\in\Z^d$},
%
%e2.14 ###
\begin{equation}\label{b-eq:def_gT_d>2}
g_T(z) = (1+|z|)^{2-d}\min\biggl\{1,\biggl(\frac{|z|}{\sqrt{T}}\biggr)^{-3}\biggr\},
\end{equation}
and for $d=2$,
%
%e2.15 ###
\begin{equation}\label{b-eq:def_gT_d=2}
g_T(z) = \left\{
\matrix{
\displaystyle \ln\biggl(\frac{1+|z|}{\sqrt{T}}\biggr)\mbox{ for
}|z|\leq\sqrt{T}\vspace*{2pt}\cr
\displaystyle \biggl(\frac{|z|}{\sqrt{T}}\biggr)^{-3}\mbox{ for }|z|>\sqrt{T}
}
\right\}.
\end{equation}
Then we have
%
%e2.16 ###
\begin{equation} \label{b-eq:dyad}
%&&
\int_{\Z^d}g_T(z) \int_{\Z^d} h_T(w)h_T(z-w)\,dw\,dz
\lesssim
\cases{
d=2,&\quad $ T$,\cr
d=3,&\quad $ \sqrt{T}$,\cr
d=4,&\quad $ \ln T$, \cr
% \ln(T\ln T),\\
d>4,&\quad $  1$.
}
\end{equation}
\end{lemma}
%

%%%%%%%%%%%%%%%%%%%%%%%%%%%%%%%%%%%%%%%%%%%%%%%%%%%%%%%%%
%%%%%%%%%%%%%%%%%%%%%%%%%%%%%%%%%%%%%%%%%%%%%%%%%%%%%%%%%

%s3 ###
\section{Proof of the main results}

Throughout this section, we let $\xi\in\R^d$ be such that\vadjust{\goodbreak}
$|\xi|=1$.

%s3.1 ###
\subsection{\texorpdfstring{Proof of Theorem \protect\ref{th1}}{Proof of Theorem 1}}

In view of \eqref{eq:err-syst}, in order to estimate $|A_T-A_\ho|$,
we need to estimate how close the modified corrector $\phi_T$
is to the original corrector $\phi$ [in terms of
$\langle|\nabla\phi_T-\nabla\phi|^2 \rangle$]. Therefore, it is
natural to introduce $\psi_T=T^2\frac{\partial\phi_T}{\partial T}$
(the prefactor $T^2$ is such that $\psi_T$ is properly renormalized
in the limit $T\uparrow\infty$ at least for large $d$).
Considering $\psi_T$ is also convenient
since for $d=2$, the corrector $\phi$ is not known to be stationary
(only its gradient is known to be stationary) so that
working with the modified correctors $\phi_T$,
which are known to be stationary, avoids technical subtleties.
In fact, we opt for a
dyadically discrete version of $\psi_T$ defined via
%
%e3.1 ###
\begin{equation}\label{eq:fo-psi}
\psi_T := T(\phi_{2T}-\phi_T).
\end{equation}
This discrete version has the technical advantage that we do not have to
think about the differentiability
of $\phi_T$ in $T$.
Moreover, its dyadic nature is in line with the dyadic decomposition of
the $T$-axis according to
%
%e3.2 ###
\begin{equation}\label{eq:dyad-Taxis}
|A_T-A_{\ho}| \le\sum_{i=0}^\infty|A_{2^iT}-A_{2^{i+1}T}|
\end{equation}
forced upon us in the case of $d=2$.
In order to get \eqref{eq:dyad-Taxis}, we used the fact that
%
%e3.3 ###
\begin{equation}\label{eq:limAT=Ahom}
\lim_{T\to\infty}A_T = A_\ho,
\end{equation}
which is proved in \cite{Gloria-Otto-09}, Proof of Theorem~1, Step~8.
We shall also use that~$\psi_T$ solves
%
%e3.4 ###
\begin{equation}\label{eq:def-psiTU}
T^{-1}\psi_{T}-\nabla^*\cdot A\nabla\psi_{T} = \tfrac{1}{2}\phi_{2T}.
\end{equation}

We split the proof in eight steps.

\step{1} Derivation of
%
%e3.5 ###
\begin{equation}\label{eq:step2-main}
|\xi\cdot(A_{2T}-A_{T})\xi| \leq T^{-2}|\langle\phi_T\psi _{T}
\rangle|+\frac{T^{-2}}{2}|\langle\phi_{2T}\psi_{T} \rangle|.
\end{equation}
Although this could be directly inferred from the spectral formula
\eqref{intro:spectral} for~$A_T$,
we give an elementary argument relying only on the corrector equation.
We recall the following consequence of (\ref{eq:app-corr})
which is proved in \cite{Gloria-Otto-09}, Proof of Theorem~1, Step~8:
%
%e3.6 ###
\begin{equation}\label{eq:VF-proba}
T^{-1}\langle\phi_T\chi \rangle+\langle(\xi+\nabla\phi_T)\cdot
A\nabla \chi \rangle = 0
\end{equation}
for every field $\chi\dvtx\Z^d\to\R$ that is jointly stationary with
$A$ and such that $\langle\chi^2 \rangle<\infty$.
From formally differentiating the definition \eqref{eq:def-AT} of
$A_T$ w.r.t. $T$ and using (\ref{eq:VF-proba})
for $\chi=\frac{\partial\phi_T}{\partial T}$, we obtain
\[
\xi\cdot\frac{\partial A_T}{\partial T}\xi
= -2T^{-1}\biggl\langle\frac{\partial\phi_T}{\partial T}\phi_T\biggr\rangle.
\]
We claim that the corresponding discrete-in-$T$ version reads
%
%e3.7 ###
\begin{equation}\label{S1.1}
\xi\cdot(A_{2T}-A_T)\xi
= -T^{-2} \bigl(\langle\psi_T\phi_T\rangle
+\tfrac{1}{2}\langle\psi_T\phi_{2T}\rangle\bigr).\vadjust{\goodbreak}
\end{equation}
Indeed, by definition of $A_T$, by expanding the square, by symmetry of
$A$, by definition of $\psi_T$,
and (\ref{eq:VF-proba}), we have
\begin{eqnarray*}
&&\xi\cdot(A_{2T}-A_T)\xi\\
&&\hspace*{3.8pt}\qquad =
\langle(\xi+\nabla\phi_{2T})\cdot A(\xi+\nabla\phi_{2T})\rangle
-\langle(\xi+\nabla\phi_{T})\cdot A(\xi+\nabla\phi_{T})\rangle
\\
&&\hspace*{3.8pt}\qquad =
\langle(\nabla\phi_{2T}-\nabla\phi_T)\cdot A(\xi+\nabla\phi
_{2T})\rangle
+\langle(\nabla\phi_{2T}-\nabla\phi_{T})\cdot A(\xi+\nabla\phi
_{T})\rangle\\
&&\qquad \stackrel{\mbox{\fontsize{8.36}{9}\selectfont\eqref{eq:fo-psi}}}{=}
T^{-1} \bigl(\langle\nabla\psi_T\cdot A(\xi+\nabla\phi_{2T})\rangle
+\langle\nabla\psi_T\cdot A(\xi+\nabla\phi_{T})\rangle\bigr)
\\
&&\qquad \stackrel{\mbox{\fontsize{8.36}{9}\selectfont(\ref{eq:VF-proba})}}{=}
-T^{-1} \bigl((2T)^{-1}\langle\psi_T\phi_{2T}\rangle
+T^{-1}\langle\psi_T\phi_{T}\rangle\bigr).
\end{eqnarray*}

In the next four steps, we focus on the first term of the
r.h.s. of \eqref{eq:step2-main}.
The second term will be dealt with the same way in Step~7.

\step{2} Proof of
%
%e3.8 ###
\begin{equation}  \label{eq:step3-main}
|\langle\phi_T\psi_{T} \rangle|
\lesssim \sum_e \biggl\langle\sup_{a(e)} \biggl|\frac{\partial\phi
_T(0)}{\partial a(e)} \biggr|^2 \biggr\rangle^{1/2} \biggl\langle\sup_{a(e)} \biggl|\frac
{\partial \psi_{T}(0)}{\partial a(e)} \biggr|^2 \biggr\rangle^{1/2},
\end{equation}
where the sum runs over the edges $e$, and proof of the representation formulas
%
%e3.9 ###
\begin{eqnarray} \label{eq:step4-main-1}
\frac{\partial\phi_T(0)}{\partial a(e)}&=&-\bigl(\xi_i+\nabla_i\phi
_T(z)\bigr)\nabla_{z_i}G_T(z,0) ,\\[-2pt]
\frac{\partial\psi_{T}(0)}{\partial a(e)}&=&-\nabla_i\psi
_{T}(z)\nabla_{z_i}G_T(z,0)
\label{eq:step4-main-2} \nonumber\\[-9pt]\\[-9pt]
&&{} -\frac{1}{2}\int_{\Z^d}G_T(0,w)\bigl(\xi_i+\nabla_i\phi
_{2T}(z)\bigr)\nabla_{z_i}G_{2T}(z,w)\,dw,\nonumber
\end{eqnarray}
where the\vadjust{\goodbreak} edge is $e=[z,z+\ee_i]$.

Due to \cite{Gloria-Otto-09}, Lemma~2.6, the functions $\phi
_T$ and $\psi_{T}$ are measurable with respect to the coefficients $a$.
Hence, \eqref{eq:step3-main} is a consequence of the covariance
estimate of Lemma \ref{lem:cov}: Since $\langle\phi_T \rangle
=\langle\psi_T \rangle=0$,
\begin{eqnarray*}
\langle\phi_T\psi_{T} \rangle&=&\bigl\langle(\phi_T-\langle\phi_T
\rangle)(\psi _{T}-\langle\psi_{T} \rangle) \bigr\rangle \\[-2pt]
&=& \cov{\phi_T}{\psi_{T}}.
\end{eqnarray*}

Formula \eqref{eq:step4-main-1} is identical to \cite
{Gloria-Otto-09}, Lemma~2.4,
(2.12).
To prove \eqref{eq:step4-main-2}, we first make use of the Green
representation formula for the solution to \eqref{eq:def-psiTU}:
%
%e3.10 ###
\begin{equation}\label{eq:step4-main-psi}
\psi_{T}(x) = \frac{1}{2}\int_{\Z^d}G_T(x,w)\phi_{2T}(w)\,dw
\end{equation}
for all $x\in\Z^d$.
Since $a(e)\mapsto\phi_{T}(\cdot;a(e))$ and $a(e)\mapsto\phi
_{2T}(\cdot;a(e))$
are continuously differentiable by \cite{Gloria-Otto-09}, Lemma~2.4,
we deduce by formula \eqref{eq:fo-psi} that $a(e)\mapsto\psi
_{T}(\cdot;a(e))$ is also continuously differentiable.
Using then the formulas \cite{Gloria-Otto-09}, Lemma~2.5, (2.15), and
\cite{Gloria-Otto-09}, Lemma~2.4, (2.12), for the derivatives of $G_T$
and $\phi_T$ with respect to $a(e)$,
and the fact that $G_T\in L^1(\Z^d)$ (see~\cite
{Gloria-Otto-09}, Corollary~2.2),
we may switch the order of the differentiation and the integration to
obtain for all $x\in\Z^d$
%
%e3.11 ###
\begin{eqnarray}\label{eq:diff-psiT}
&&\frac{\partial\psi_{T}(x)}{\partial a(e)}
= \frac{1}{2}\int_{\Z^d}\frac{\partial G_T(x,w)}{\partial
a(e)}\phi_{2T}(w)\,dw\nonumber \\[-2pt]
&&\hphantom{\frac{\partial\psi_{T}(x)}{\partial a(e)}
=}{}+ \frac{1}{2}\int_{\Z^d}G_T(x,w)\frac{\partial
\phi_{2T}(w)}{\partial a(e)}\,dw \nonumber\\[-2pt]
&&\hspace*{7.4pt}\stackrel{\mbox{\fontsize{8.36}{9}\selectfont\cite{Gloria-Otto-09}}, (2.12)\ \mathrm{and}\
(2.15)}{=}-\frac
{1}{2}\int_{\Z^d}\nabla_{z_i}G_T(x,z)\nabla_{z_i}G_T(z,w) \phi
_{2T}(w)\,dw \\[-2pt]
&&\hspace*{46pt} {} -\frac{1}{2}\int_{\Z^d}G_T(x,w)\bigl(\xi_i+\nabla_i \phi
_{2T}(z)\bigr)\nabla_{z_i}G_{2T}(z,w)\,dw \nonumber\\[-2pt]
&&\hspace*{33.5pt}\stackrel{\mbox{\fontsize{8.36}{9}\selectfont\eqref{eq:step4-main-psi}}}{=} -\nabla
_{z_i}G_T(x,z)\nabla_i \psi_{T}(z)\nonumber \\[-2pt]
&&\hspace*{46pt} {}- \frac{1}{2}\int_{\Z
^d}G_T(x,w)\bigl(\xi_i+\nabla_{i}\phi_{2T}(z)\bigr)\nabla_{z_i}
G_{2T}(z,w)\,dw,\nonumber
\end{eqnarray}
which is \eqref{eq:step4-main-2} taking $x=0$.

From now on in the proof, we let $g_T$ be defined as in
Lemma~\ref{b-lem:dyad} (i.e., $g_T$ decays as the Green function $G_T$).

\step{3} In this step, we shall prove that
%
%e3.12 ###
\begin{eqnarray}
|\langle\phi_T\psi_{T} \rangle|&\lesssim&\mathcal{L} + \mathcal
{N}\label
{eq:step5-main},
\end{eqnarray}
where
\begin{eqnarray}\label{eq:L-def}
\qquad \mathcal{L}&:=& \int_{\Z^d} \bigl\langle\bigl(1+|\nabla\phi
_T(z)|^2\bigr)|\nabla _z G_T(z,0)|^2
\bigr\rangle^{1/2}\nonumber\\[-9pt]\\[-9pt]
&&\hphantom{\int_{\Z^d}}{}\times \biggl\langle\biggl(\int_{\Z^d}g_T(w)\bigl(1+|\nabla\phi _{2T}(z)|\bigr)|\nabla
_zG_{2T}(z,w)|\,dw \biggr)^2 \biggr\rangle^{1/2}\,dz,\nonumber
\end{eqnarray}
and $\mathcal{N}=\mathcal{N}_1+\mathcal{N}_2$,
\begin{eqnarray}\label{eq:N-def}
\mathcal{N}_1&:=& \int_{\Z^d} \bigl\langle\bigl(1+|\nabla\phi
_T(z)|^2\bigr)|\nabla_z G_T(z,0)|^2 \bigr\rangle^{1/2}\nonumber\\[-8pt]\\[-8pt]
&&\hphantom{\int_{\Z^d}}
{}\times\langle|\nabla\psi_{T}(z)|^2|\nabla_zG_T(z,0)|^2 \rangle^{1/2}\,dz,\nonumber
\\
\label{eq:N2-def}
\mathcal{N}_2&:=& \mu_d(T)\int_{\Z^d} \bigl\langle\bigl(1+|\nabla\phi
_T(z)|^2\bigr)|\nabla_z G_T(z,0)|^2 \bigr\rangle^{1/2} \nonumber\\[-8pt]\\[-8pt]
&& \hphantom{\mu_d(T)\int_{\Z^d}}{}\times\bigl\langle\bigl(1+|\nabla\phi_{2T}(z)|^2\bigr)|\nabla _zG_T(z,0)|^2 \bigr\rangle
^{1/2}\,dz\nonumber
\end{eqnarray}
with
\[
\mu_d(T) := \cases{
d=2,&\quad $\ln T$, \cr
d>2,&\quad 1.
}
\]
The term $\mathcal{L}$ is a linear error: It is of the same type as
for the analysis in the limit of vanishing ellipticity contrast
(see \cite{Gloria-Otto-09}, the Appendix).
On the contrary, the term $\mathcal{N}$ is nonlinear and does not
appear in the limit of vanishing ellipticity contrast.
As we shall prove, it is of lower order.
The terms $\mathcal{L}$ and $\mathcal{N}_1$ in estimate \eqref
{eq:step5-main} would be direct consequences of \eqref{eq:step3-main},
and \eqref{eq:step4-main-1} and \eqref{eq:step4-main-2},
disregarding the suprema in $a(e)$ in \eqref{eq:step3-main}. Taking
the suprema in $a(e)$ into account actually brings the second nonlinear
term $\mathcal N_2$,
which turns out to be of lower order than $\mathcal{N}_1$.

According to \cite{Gloria-Otto-09}, Lemma~2.4, (2.13), we
have for \eqref{eq:step4-main-1}
%
%e3.13 ###
\begin{eqnarray}\label{eq:sup-1}
\sup_{a(e)} \biggl|\frac{\partial\phi_T(0)}{\partial a(e)} \biggr|&\lesssim
&\bigl(1+|\nabla_i\phi_T(z)|\bigr)|\nabla_{z}G_T(z,0)|.
\end{eqnarray}
It remains to deal with \eqref{eq:step4-main-2}.
Using the pointwise decay of $G_T$ in Lemma~\ref
{lem:ptwise-Green-decay} combined
with the susceptibility estimates \cite{Gloria-Otto-09}, Lemma~2.4,
(2.14), and \cite{Gloria-Otto-09}, Lemma~2.5, (2.16),
of $\nabla\phi_T$ and $\nabla G_T$
w.r.t. $a(e)$, we obtain
\begin{eqnarray}\label{eq:sup-2}
&&\sup_{a(e)} \biggl|\frac{1}{2}\int_{\Z^d}G_T(0,w)\bigl(\xi_i+\nabla
_i\phi_{2T}(z)\bigr)\nabla_{z_i}G_{2T}(z,w)\,dw \biggr|\nonumber\\[-8pt]\\[-8pt]
&&\qquad \lesssim \int_{\Z^d}g_T(w)\bigl(1+|\nabla_i\phi_{2T}(z)|\bigr)|\nabla
_{z}G_{2T}(z,w)|\,dw,\nonumber
\end{eqnarray}
which together with \eqref{eq:sup-1} gives the linear term ${\mathcal L}$.

To treat the first term of the r.h.s. of \eqref{eq:step4-main-2},
we need to deal with the supremum of $|\nabla_i \psi_{T}(z)|$ over $a(e)$.
We appeal to \eqref{eq:diff-psiT} that we rewrite in the form
\begin{eqnarray*}
\frac{\partial\psi_{T}(x)}{\partial a(e)}
&=&-\nabla_i \psi_T(z) G_T(x,e)\\
&&{}-\frac{1}{2}\bigl(\xi_i+\nabla_i \phi
_{2T}(z)\bigr)\int_{\Z^d}G_T(x,w) G_T(e,w)\, dw ,
\end{eqnarray*}
where $G_T(x,e):=G_T(x,z+\ee_i)-G_T(x,z)$ and $G_T(e,w):=G_T(z+\ee
_i,w)-G_T(z,w)$.
Hence,
\begin{eqnarray}\label{eq:step5-estimDPsi}
\frac{\partial\nabla_i \psi_{T}(z)}{\partial a(e)}&=& -\nabla_i
\psi_T(z) G_T(e,e)\nonumber\\[-8pt]\\[-8pt]
&&{}-\frac{1}{2} \bigl(\xi_i+\nabla_i \phi_{2T}(z)\bigr) \int
_{\Z^d}G_T(e,w) G_{2T}(e,w)\, dw,\nonumber
\end{eqnarray}
where $G_T(e,e):=G_T(z+\ee_i,z+\ee_i)+G_T(z,z)-G_T(z+\ee
_i,z)-G_T(z,z+\ee_i)$.
On the one hand, the uniform bound \cite{Gloria-Otto-09},
Corollary~2.3, on $\nabla G_T$ yields
$|G_T(e,e)|\lesssim1$.
On the other hand, as we shall argue, the integrability of~$\nabla G_T$
and $\nabla G_{2T}$ from \cite{Gloria-Otto-09},\vadjust{\goodbreak} Lemma~2.9 (combined
with the uniform bound~\cite{Gloria-Otto-09}, Corollary~2.3, on
gradients) implies
%
%e3.14 ###
\begin{equation} \label{eq:th1-step5-grad2}
\int_{\Z^d}G_T(e,w) G_{2T}(e,w)\, dw \lesssim \mu_d(T) = \cases{
d=2,&\quad $ \ln T $,\cr
d>2,&\quad 1.
}
\end{equation}
Hence, if we regard \eqref{eq:step5-estimDPsi} as an ordinary
differential equation for $\nabla_i\psi_T(z)$ in the variable $a(e)$,
we obtain
%
%e3.15 ###
\begin{eqnarray}\label{eq:sup-psi}
\sup_{a(e)} | \nabla_i \psi_{T}(z) |&\lesssim& |\nabla_i \psi
_T(z)| + \mu_d(T)\bigl(1+|\nabla_i \phi_{2T}(z)|\bigr)
\end{eqnarray}
since $a(e)$ lies in a bounded domain $[\alpha,\beta]$, and $\sup
_{a(e)}|\nabla_i \phi_{2T}(z)|\lesssim1+|\nabla_i \phi_{2T}(z)|$
according to \cite{Gloria-Otto-09}, Lemma 2.4, (2.14), with $2T$
instead of $T$.
Note that \eqref{eq:sup-1}, \eqref{eq:sup-psi} and $\sup
_{a(e)}|\nabla_{z_i}G_T(z,0)|\lesssim|\nabla_{z_i}G_T(z,0)|$ give
the nonlinear terms $\mathcal{N}_1$ and $\mathcal{N}_2$.

We now give the argument for \eqref{eq:th1-step5-grad2}. We
first use the Cauchy--Schwarz inequality
\begin{eqnarray*}
&&\int_{\Z^d}G_T(e,w) G_{2T}(e,w) \,dw\\
&&\qquad \leq \biggl(\int_{\Z^d}G_T(e,w)^2dw
\biggr)^{1/2} \biggl( \int_{\Z^d}G_{2T}(e,w)^2\, dw \biggr)^{1/2} \\
&&\qquad \leq \biggl(\int_{\Z^d}|\nabla_z G_T(z,w)|^2\,dw \biggr)^{1/2} \biggl( \int_{\Z
^d}|\nabla_z G_{2T}(z,w)|^2 \,dw \biggr)^{1/2}
\end{eqnarray*}
and then make a decomposition of $\Z^d$ into the ball
of radius $R\sim1$, and dyadic annuli $\{w\dvtx 2^i R < |z-w|\leq2^{i+1}R\}
$ for $i\in\N$.
On the ball of radius~$R$, we use the uniform estimate of \cite
{Gloria-Otto-09}, Corollary 2.3, on $\nabla G_T$,
whereas on the dyadic annuli we appeal to the decay estimate in \cite
{Gloria-Otto-09}, Lemma 2.9, for the gradient
of the Green function, which requires $R$ to be sufficiently large
although still of order $1$.
Both terms in the r.h.s. scale the same way and we only treat the first one:
\begin{eqnarray*}
&&\int_{\Z^d}|\nabla_zG_T(z,w)|^2\,dw\\
&&\qquad =\int_{|z-w|\leq R}|\nabla_zG_T(z,w)|^2\,dw+\sum_{i=0}^\infty\int
_{2^iR<|z-w|\leq2^{i+1}R}|\nabla_z G_T(z,w)|^2\,dw\\
&&\qquad \lesssim 1+\sum_{i=1}^\infty(2^i)^{d+2(1-d)}\min\bigl\{1,\sqrt
{T}(2^{i}R)^{-1}\bigr\}^2\\
&&\qquad \lesssim\mu_d(T),
\end{eqnarray*}
using \cite{Gloria-Otto-09}, Corollary 2.3, and \cite
{Gloria-Otto-09}, Lemma 2.9, for $k=2$, respectively. This concludes
Step 3.

\step{4} Suboptimal estimate of the nonlinear term $\mathcal{N}$:
%
%e3.17 ###
%e3.16 ###
\begin{eqnarray}\label{eq:step6-main1}
\mathcal{N}_1&\lesssim&
\langle\nabla\psi_{T}\cdot A\nabla\psi_{T} \rangle ^{1/2} \cases{
d=2,&\quad $\sqrt{T} \ln^q T$, \vspace*{-1pt}\cr
d=3,&\quad $\ln T$, \vspace*{-1pt}\cr
d>3,&\quad 1,
}
\\[-1pt]
\mathcal{N}_2&\lesssim& \mu_d(T)^q,\label{eq:step6-main2}
\end{eqnarray}
where $q$ is a generic exponent which only depends on $\alpha,\beta$.
We first deal with $\mathcal{N}_1$, and begin with the second factor
of the r.h.s. of \eqref{eq:N-def}.
The pointwise estimate \eqref{eq:def_gT_d>2} of Lemma~\ref
{lem:ptwise-Green-decay} for $d>2$
on the Green function gives the \emph{suboptimal pointwise} estimate
on the gradient of the Green function
%
%e3.18 ###
\begin{equation}\label{eq:nabla-par-G}
|\nabla G_T(z,0)|\leq G_T(z,0)+\sum_{i=1}^d G_T(z+\ee_i,0)
\lesssim(1+|z|)^{2-d}.
\end{equation}
This estimate coincides for $d=2$ with the uniform bound of \cite
{Gloria-Otto-09}, Corollary~2.3.
The coercivity of $A$ thus yields
\begin{eqnarray*}
&&\langle|\nabla G_T(z,0)|^2|\nabla\psi_{T}(z)|^2 \rangle
^{1/2} \\[-2pt]
&&\qquad \lesssim(1+|z|)^{2-d}\langle\nabla\psi_{T}(z)\cdot A(z)\nabla
\psi_{T}(z) \rangle^{1/2}\\[-2pt]
&&\qquad =(1+|z|)^{2-d}\langle\nabla\psi_{T}\cdot A\nabla\psi_{T} \rangle^{1/2}
\end{eqnarray*}
by joint stationarity of $\nabla\psi_T$ and $A$.
Hence, \eqref{eq:N-def} turns into
\[
\mathcal{N}_1 \lesssim\langle\nabla\psi_{T}\cdot A\nabla\psi
_{T} \rangle^{1/2}\int_{\Z^d}(1+|z|)^{2-d}
\bigl\langle\bigl(1+|\nabla\phi_T(z)|^2\bigr)|\nabla G_T(z,0)|^2 \bigr\rangle^{1/2}\,dz.
\]
We then let $p>2$ be a Meyers' exponent as in \cite
{Gloria-Otto-09}, Lemma~2.9 and use H\"{o}lder's inequality
in probability with exponents $(p/(p-2),p/2)$, the stationarity of
$\nabla\phi_T$, the fact that the gradient of $\phi_T$ is estimated
by $\phi_T$
as in \eqref{eq:nabla-par-G}, and the bounds on the stochastic moments
of $\phi_T$ in \cite{Gloria-Otto-09}, Proposition~1,\looseness=-1
%
%e3.19 ###
\begin{eqnarray}\label{eq:modif-N}
\quad \mathcal{N}_1&\lesssim& \langle\nabla\psi_{T}\cdot A\nabla\psi
_{T} \rangle^{1/2}\nonumber\\[-2pt]
&&{}\times \int_{\Z^d}(1+|z|)^{2-d}
\bigl\langle1+|\nabla\phi_T(z)|^{2p/(p-2)} \bigr\rangle^{(p-2)/(2p)}\langle
|\nabla G_T(z,0)|^p \rangle^{1/p}\,dz \nonumber\\[-2pt]
&=& \langle\nabla\psi_{T}\cdot A\nabla\psi_{T} \rangle
^{1/2}\bigl\langle1+|\nabla\phi_T|^{2p/(p-2)} \bigr\rangle^{(p-2)/(2p)}
\nonumber\\[-2pt]
&&{}\times \int_{\Z^d}(1+|z|)^{2-d}
\langle|\nabla G_T(z,0)|^p \rangle^{1/p}\,dz \\[-2pt]
&\lesssim&\langle\nabla\psi_{T}\cdot A\nabla\psi_{T} \rangle
^{1/2}\bigl\langle1+|\phi_T|^{2p/(p-2)} \bigr\rangle^{(p-2)/(2p)} \nonumber\\[-2pt]
&&{}\times \int_{\Z
^d}(1+|z|)^{2-d}
\langle|\nabla G_T(z,0)|^p \rangle^{1/p}\,dz \nonumber\\[-2pt]
&\lesssim& \mu_d(T)^q\langle\nabla\psi_{T}\cdot A\nabla\psi _{T}
\rangle^{1/2} \int_{\Z^d}(1+|z|)^{2-d}
\langle|\nabla G_T(z,0)|^p \rangle^{1/p}\,dz,\nonumber
\end{eqnarray}\looseness=0
for some generic $q$ depending only on $\alpha,\beta$.
H\"{o}lder's inequality with exponents\vadjust{\goodbreak} $(p, p/(p-1))$ in $\Z^d$,
combined with the same dyadic decomposition of $\Z^d$
as for the proof of \eqref{eq:th1-step5-grad2}
(and the uniform bound on $\nabla G_T$ from \cite
{Gloria-Otto-09}, Corollary~2.3) yields
\begin{eqnarray*}
&&\int_{\Z^d}(1+|z|)^{2-d} \langle|\nabla G_T(z,0)|^p\rangle^{1/p} \,dz \\
&&\qquad \lesssim
1+\sum_{i=0}^\infty
\biggl(\biggl\langle\int_{2^iR< |z|\leq2^{i+1}R} |\nabla G_T(z,0)|^p\,dz \biggr\rangle \biggr)^{1/p}
\\
&&\hphantom{1+\sum_{i=0}^\infty}\qquad\quad
{}\times\biggl(\int_{2^iR< |z|\leq2^{i+1}R} (1+|z|)^{(2-d)p/(p-1)}\,dz \biggr)^{(p-1)/p}.
\end{eqnarray*}
Using the optimal decay of $\nabla G_T$ on dyadic annuli in $L^p$ norm
from \cite{Gloria-Otto-09}, Lemma~2.9, with $k=2p$, this turns into
\begin{eqnarray*}
&&\int_{\Z^d}(1+|z|)^{2-d} \langle|\nabla G_T(z,0)|^p\rangle^{1/p}\,dz\\
&&\qquad \lesssim
1+\sum_{i=0}^\infty\biggl( (2^iR)^{d}(2^iR)^{(1-d)p}\min\biggl\{1,\frac{\sqrt
{T}}{2^iR}\biggr\}^{2p} \biggr)^{1/p} \\
&&\qquad \quad
\hphantom{1+\sum_{i=0}^\infty}{}\times\bigl( (2^iR)^{d}(2^iR)^{(2-d)p/(p-1)}
\bigr)^{(p-1)/p}\\
&&\qquad =
1+\sum_{i=0}^\infty(2^iR)^{3-d}\min\biggl\{1,\frac{\sqrt{T}}{2^iR}\biggr\}^2.
\end{eqnarray*}
Recalling that $R\sim1$, this implies
\begin{eqnarray*}
\int_{\Z^d}(1+|z|)^{2-d} \langle|\nabla G_T(z,0)|^p \rangle^{1/p}
\,dz&\lesssim&
\cases{
d=2,&\quad $ \sqrt{T}$, \cr
d=3,&\quad $ \ln T$, \cr
d>3,&\quad  1.
}
\end{eqnarray*}
Combined with \eqref{eq:modif-N} it proves \eqref{eq:step6-main1}.

% Combined with \eqref{eq:modif-N} and Young's inequality, it then
%yields
% %
% \begin{equation}\label{eq:step6-main-part1}
% \mathcal{N}_1-\frac{1}{2}\expec{\nabla\psi_{T}\cdot A\nabla
% \{
% \begin{array}{lcl}
% d=2&:&{T} (\ln T)^q, \\
% d=3&:&\ln^2 T, \\
% d>3&:&1.
% \end{array}
%
% \end{equation}
% %

We now turn to $\mathcal{N}_2$.
Proceeding as above to deal with the terms $\nabla\phi_T$ and~$\nabla
\phi_{2T}$ in $\mathcal{N}_2$, we obtain as desired
\begin{eqnarray*}
\mathcal{N}_2 &\lesssim& \mu_d(T)\mu_d(T)^{2q}\int_{\Z^d}\langle
|\nabla_zG_T(z,0)|^p \rangle^{2/p}\,dz \nonumber\\
&\lesssim& \mu_d(T)^{2q+2},
\end{eqnarray*}
using the same dyadic decomposition of $\Z^d$ as for the proof of
\eqref{eq:th1-step5-grad2} together with the higher integrability of
gradients of
\cite{Gloria-Otto-09}, Lemma~2.9 and \cite{Gloria-Otto-09}, Corollary~2.3.

\step{5} Estimate of the linear term $\mathcal{L}$:
%
%e3.20 ###
\begin{equation}\label{eq:step7-main}
\mathcal{L} \lesssim
\cases{
d=2,&\quad $T \ln^qT $,\cr
d=3,&\quad $\sqrt{T} $,\cr
d=4,&\quad $\ln T $,\cr
d>4,&\quad 1.
}
\end{equation}
We first treat the second factor of \eqref{eq:L-def}.
We proceed as in Step~4 to deal with the expectation of the corrector
term, and
let $p>2$ be a Meyers' exponent as in \cite{Gloria-Otto-09}, Lemma~2.9.
We obtain by H\"{o}lder's inequality in probability with exponents
$(p/(p-2),p,p)$ and the bounds on the stochastic
moments of $\phi_T$ from \cite{Gloria-Otto-09}, Proposition~1:
\begin{eqnarray*}
&&\biggl\langle\biggl(\int_{\Z^d}g_T(w)\bigl(1+|\nabla\phi
_{2T}(z)|\bigr)|\nabla_{z_i}G_{2T}(z,w)|\,dw \biggr)^2 \biggr\rangle\\
&&\qquad =\int_{\Z^d}\int_{\Z^d}g_T(w)g_T(w')\\
&&\quad\qquad \hphantom{\int_{\Z^d}\int_{\Z^d}}
{}\times\bigl\langle\bigl(1+|\nabla\phi
_{2T}(z)|\bigr)^2|\nabla_{z_i}G_{2T}(z,w)||\nabla_{z_i}G_{2T}(z,w')|
\bigr\rangle \,dw\,dw'
\\
&&\qquad \lesssim \bigl(1+\bigl\langle|\phi_{2T}|^{2p/(p-2)} \bigr\rangle^{(p-2)/p} \bigr) \\
&&\qquad \quad {}\times \int_{\Z^d}\int_{\Z^d}g_T(w)g_T(w')\langle|\nabla
_{z_i}G_{2T}(z,w)|^{p} \rangle^{1/p}\\
&&\qquad \quad \hphantom{{}\times \int_{\Z^d}\int_{\Z^d}}
{}\times\langle|\nabla
_{z_i}G_{2T}(z,w')|^{p} \rangle^{1/p}\,dw\,dw' \\
&&\qquad \lesssim \mu_d(T)^q \biggl(\int_{\Z^d}g_T(w)\langle|\nabla_z
G_{2T}(z,w)|^p \rangle^{1/p}\,dw \biggr)^2.
\end{eqnarray*}
We thus have
\begin{eqnarray*}
\mathcal{L}
&\lesssim&\mu_d(T)^q \int_{\Z^d}\bigl\langle\bigl(1+|\nabla\phi
_T(z)|^2\bigr)|\nabla G_T(z,0)|^2 \bigr\rangle^{1/2}\\
&&\hphantom{\mu_d(T)^q \int_{\Z^d}}{}\times\int_{\Z^d}g_T(w)\langle
|\nabla _z G_{2T}(z,w)|^p \rangle^{1/p}\,dw\,dz.
\end{eqnarray*}
Appealing once more to H\"{o}lder's inequality in probability with
exponents $(p/(p-2),p/2)$ and to
\cite{Gloria-Otto-09}, Proposition~1, this turns into
\begin{eqnarray*}
\mathcal{L}
&\lesssim&\mu_d(T)^{2q}\int_{\Z^d}g_T(w)\int_{\Z^d}\langle
|\nabla_z G_{2T}(z,w)|^p \rangle^{1/p}\langle|\nabla_z G_T(z,0)|^p
\rangle^{1/p}\,dz\,dw \\
&=&\mu_d(T)^{2q}\int_{\Z^d}g_T(w)\int_{\Z^d}h_{2T}(z-w)h_T(z)\,dz\,dw,
\end{eqnarray*}
where, by stationarity, we have set
\begin{eqnarray*}
h_T(w) &=& \langle|\nabla_w G_T(w,0)|^p \rangle^{1/p},\\
h_{2T}(w)&= &\langle|\nabla_w G_{2T}(w,0)|^p \rangle^{1/p}.
\end{eqnarray*}
By the optimal decay estimate of $\nabla G_T$ on dyadic annuli from
\cite{Gloria-Otto-09}, Lemma~2.9 (and by the uniform bounds
on $\nabla G_T$ from \cite{Gloria-Otto-09}, Corollary~2.3),
and by definition of $g_T$, we are in position to apply Lemma~\ref{b-lem:dyad}.
Estimate \eqref{eq:step7-main} is thus proved.

\step{6} Proof of
%
%e3.21 ###
\begin{equation}\label{eq:step1bis-main}
\langle\nabla\psi_{T}\cdot A \nabla\psi_{T} \rangle \leq|\langle
\phi _T\psi_{T} \rangle|.
\end{equation}
Using \eqref{eq:fo-psi}, we rewrite \eqref{eq:def-psiTU} as
\begin{eqnarray}\label{eq:new-psiTU}
(2T)^{-1}\psi_{T}-\nabla^*\cdot A\nabla\psi_{T}&=&\tfrac{1}{2}\phi
_{2T}-(2T)^{-1}\psi_{T}\nonumber\\[-8pt]\\[-8pt]
&=&\tfrac{1}{2} \phi_T.\nonumber
\end{eqnarray}
We now multiply \eqref{eq:new-psiTU} by $\psi_{T}$:
\begin{eqnarray*}
(2T)^{-1}\psi_{T}^2-(\nabla^*\cdot A\nabla\psi_{T})\psi
_{T}&=&\tfrac{1}{2}\phi_T\psi_{T}.
\end{eqnarray*}
By integration by parts and joint stationarity of $\psi_{T}$, $\nabla
\psi_{T}$ and $A$ (see \cite{Gloria-Otto-09}, Proof of Theorem~1,
Step~8, for details),
this turns into
\begin{eqnarray*}
(2T)^{-1}\langle\psi_{T}^2 \rangle+\langle\nabla\psi_{T}\cdot A
\nabla \psi_{T} \rangle&=&\tfrac{1}{2} \langle\phi_T\psi_{T}
\rangle.
\end{eqnarray*}
We then conclude by the nonnegativity of the first term.

\step{7} Proof of
%
%e3.22 ###
\begin{equation}\label{eq:step9-main1}
|\langle\phi_T\psi_{T} \rangle| \lesssim
\cases{
d=2,&\quad $T\ln^qT$, \cr
d=3,&\quad $\sqrt{T}$, \cr
d=4,&\quad $\ln T $,\cr
d>4,&\quad 1,
}
\end{equation}
and
%
%e3.23 ###
\begin{equation}\label{eq:step9-main2}
|\langle\phi_{2T}\psi_{T} \rangle| \lesssim
\cases{
d=2,&\quad $T \ln^qT$, \cr
d=3,&\quad $\sqrt{T}$, \cr
d=4,&\quad $\ln T $,\cr
d>4,&\quad 1.
}
\end{equation}
From Steps~3, 4 and 5, and Young's inequality, we deduce that
\[
|\langle\phi_T\psi_{T} \rangle|-\frac{1}{2}\langle\nabla\psi
_{T}\cdot A \nabla\psi_{T} \rangle \lesssim
\cases{
d=2,&\quad $T\ln^qT$, \cr
d=3,&\quad $\sqrt{T}$, \cr
d=4,&\quad $\ln T $,\cr
d>4,&\quad 1.
}
\]
Combined with Step~6, this shows \eqref{eq:step9-main1}.

For \eqref{eq:step9-main2}, we proceed exactly as for \eqref
{eq:step9-main1} in Steps~2--6.
In particular, with obvious notation, we have
\[
|\langle\phi_{2T}\psi_{T} \rangle| \lesssim\mathcal{N}'+\mathcal{L}',
\]
where
\[
\mathcal{N}'-\frac{1}{2}\langle\nabla\psi_{T}\cdot A \nabla\psi
_{T} \rangle \lesssim\cases{
d=2,&\quad $T\ln^qT$, \cr
d=3,&\quad $\ln^2T$, \cr
d>3,&\quad 1,
}
\]
and
\[
\mathcal{L}' \lesssim
\cases{
d=2,&\quad $T\ln^qT$, \cr
d=3,&\quad $\sqrt{T}$, \cr
d=4,&\quad $\ln T$,\cr
d>4,&\quad 1.
}
\]
We then conclude as above.

\step{8} Proof of \eqref{eq:estim-err-hom}.

Steps~1 and 7 yield
\begin{eqnarray}\label{eq:step10}
|\xi\cdot(A_T-A_{2T})\xi| &\leq& T^{-2}|\langle\phi_T\psi _{T}
\rangle|+(2T^2)^{-1}|\langle\phi_{2T} \psi_{T} \rangle|
\nonumber\\[-8pt]\\[-8pt]
&\lesssim& T^{-2} \cases{
d=2,&\quad $ T\ln^q T $, \cr
d=3,&\quad $\sqrt{T} $,\cr
d=4,&\quad $\ln T $,\cr
d>4,&\quad 1.
}\nonumber
\end{eqnarray}
We finally appeal to the dyadic decomposition of the $T$-axis \eqref
{eq:dyad-Taxis},
which, combined with \eqref{eq:step10}, turns into
\begin{eqnarray*}
|\xi\cdot(A_T-A_\ho)\xi| &\lesssim& \sum_{i=1}^\infty\cases{
d=2,&\quad $(2^iT)^{-1}\ln^q (2^iT)$, \cr
d=3,&\quad $(2^iT)^{-3/2}$,\cr
d=4,&\quad $ (2^iT)^{-2}\ln(2^iT)$,\cr
d>4,&\quad $(2^iT)^{-2}$,
}
\\
&\lesssim& \cases{
d=2,&\quad $ T^{-1}\ln^qT $, \cr
d=3,&\quad $T^{-3/2} $\cr
d=4,&\quad $T^{-2}\ln T $,\cr
d>4,&\quad $T^{-2} $.
}
\end{eqnarray*}
This concludes the proof of the theorem.

%s3.2 ###
\subsection{\texorpdfstring{Proof of Corollary \protect\ref{prop1}}{Proof of Corollary 1}}

By Steps~6 and 7 in the proof of Theorem~\ref{th1} and by the definition
\eqref{eq:fo-psi} of $\psi_T$, we learn that
\begin{eqnarray*}
&&\langle|\nabla\phi_{2T}-\nabla\phi_T|^2 \rangle \stackrel
{\mbox{\fontsize{8.36}{9}\selectfont\eqref
{eq:fo-psi}}}{=} T^{-2} \langle|\nabla\psi_{T}|^2 \rangle \\
&&\hspace*{61pt}\stackrel{\mbox{\fontsize{8.36}{9}\selectfont\eqref{eq:step1bis-main} and  \eqref
{eq:step9-main1}}}{\lesssim}
\cases{
d=3,&\quad $T^{-3/2}$, \cr
d=4,&\quad $T^{-2}\ln T $,\cr
d>4,&\quad $T^{-2}$.
}
\end{eqnarray*}
In particular, $\nabla\phi_T$ is a Cauchy sequence in $L^2$ in probability.
Hence, $\nabla\phi_T$ converges in $L^2$ to its weak limit $\nabla
\phi$, and
by a dyadic decomposition of the $T$-axis the above estimate yields
\[
\langle|\nabla\phi_T-\nabla\phi|^2 \rangle \lesssim\cases{
d=3,&\quad $T^{-3/2}$, \cr
d=4,&\quad $T^{-2}\ln T $,\cr
d>4,&\quad $T^{-2}$,
}
\]
which gives the second term of\vadjust{\goodbreak} the l.h.s. of \eqref{eq:prop1}.

Likewise, from Step~7 in the proof of Theorem~\ref{th1}, we
learn that
\begin{eqnarray*}
&&\langle(\phi_{2T}-\phi_T)^2 \rangle \stackrel{\mbox{\fontsize{8.36}{9}\selectfont\eqref{eq:fo-psi}}}{=}
T^{-1} \langle(\phi_{2T}-\phi_T)\psi_{T} \rangle \\
&&\hspace*{69pt}\leq T^{-1} (\langle|\phi_{2T}\psi_{T}| \rangle+\langle|\phi
_{T}\psi _{T}| \rangle )\\
&&\hspace*{46pt}\stackrel{\mbox{\fontsize{8.36}{9}\selectfont\eqref{eq:step9-main1} and \eqref
{eq:step9-main2}}}{\lesssim
}
\cases{
d=3,&\quad $T^{-1/2}$, \cr
d=4,&\quad $T^{-1}\ln T $,\cr
d>4,&\quad $T^{-1}$,
}
\end{eqnarray*}
so that $\phi_T$ is a Cauchy sequence in $L^2$ in probability and
$\phi_T$
converges in $L^2$ to its weak limit $\tilde\phi$ provided by \cite
{Gloria-Otto-09}, Corollary~1.
In particular, by a dyadic decomposition of the $T$-axis the above
estimate yields
\[
\langle(\phi_T- \tilde\phi)^2 \rangle = \cases{
d=3,&\quad $T^{-1/2}$, \cr
d=4,&\quad $T^{-1}\ln T $,\cr
d>4,&\quad $T^{-1}$,
}
\]
which is the first term of the l.h.s. of \eqref{eq:prop1}.
This concludes the proof of the corollary.

%s4 ###
\section{Proof of the auxiliary lemmas}

%s4.1 ###
\subsection{\texorpdfstring{Proof of Lemma \protect\ref{lem:cov}}{Proof of Lemma 3}}

Without loss of generality we may assume
%
%e4.1 ###
\begin{equation} \label{5.1}
\sum_{i=1}^\infty\biggl\langle\sup_{a_i} \biggl|\frac{\partial X}{\partial
a_i} \biggr|^2 \biggr\rangle, \sum_{i=1}^\infty\biggl\langle\sup_{a_i} \biggl|\frac
{\partial Y}{\partial a_i} \biggr|^2 \biggr\rangle< \infty.
\end{equation}
Let $Z_n$ denote the expected value of $Z$ conditioned on
$a_1,\ldots,a_n$, that is
\[
Z_n(a_1,\ldots,a_n) := \langle Z | a_1,\ldots,a_n \rangle.
\]
From \cite{Gloria-Otto-09}, (5.2) and (5.3), in the proof of \cite
{Gloria-Otto-09}, Lemma~2.3, we learn that
\[
\lim_{n\uparrow\infty} \langle(Z-Z_n)^2 \rangle = 0
\]
for $Z=X,Z_n=X_n$ and $Z=Y,Z_n=Y_n$, respectively,
so that, by the Cauchy--Schwarz inequality in probability,
\begin{eqnarray*}
\lim_{n\uparrow\infty} \langle X_n \rangle &=&\langle X \rangle, \\
\lim_{n\uparrow\infty} \langle Y_n \rangle &=&\langle Y \rangle, \\
\lim_{n\uparrow\infty} \langle X_nY_n \rangle &=&\langle XY \rangle.
\end{eqnarray*}
Hence,
%
%e4.2 ###
\begin{eqnarray}\label{eq:lim-covar}
\lim_{n\uparrow\infty} \cov{X_n}{Y_n} & =& \lim_{n\uparrow\infty
} (\langle X_nY_n \rangle-\langle X_n \rangle \langle Y_n \rangle)
\nonumber\\
&=& \langle XY \rangle-\langle X \rangle \langle Y \rangle
\\
&=& \cov{X}{Y}.\nonumber
\end{eqnarray}
Note also that
%
%e4.3 ###
\begin{equation}\label{5.2}
\cov{X_n}{Y_n} = \sum_{i=1}^n (\langle X_iY_i \rangle-\langle
X_{i-1}Y_{i-1} \rangle )
\end{equation}
with the notation $X_0=\langle X \rangle$ and $Y_0=\langle Y \rangle
$, so that
$\langle X_n \rangle=X_0$ and $\langle Y_n \rangle=Y_0$.
Inequality (\ref{eq:covar-estim}) then follows from \eqref{5.1},
(\ref{eq:lim-covar}), \eqref{5.2}, and
%
%e4.4 ###
\begin{eqnarray}\label{eq:pr-covar}
\langle X_iY_i\rangle-\langle X_{i-1}Y_{i-1}\rangle
&\lesssim&\biggl\langle\sup_{a_i} \biggl|\frac{\partial X}{\partial a_i} \biggr|^2
\biggr\rangle^{1/2}\biggl\langle\sup_{a_i} \biggl|\frac{\partial Y}{\partial a_i} \biggr|^2
\biggr\rangle^{1/2},
\end{eqnarray}
that we prove now.
By our assumption that $\{a_i\}_{i\in\mathbb{N}}$
are i.i.d., we have
\begin{eqnarray*}
&&Z_{i-1}(a_1,\ldots,a_{i-1})
=
\int Z_i(a_1,\ldots,a_{i-1},a_i'') \beta(da_i''),\\
&&\langle X_i(a_1,\ldots,a_i)Y_i(a_1,\ldots,a_i) \rangle\\
&&\qquad =
\biggl\langle\int X_i(a_1,\ldots,a_{i-1},a_i')Y_i(a_1,\ldots
,a_{i-1},a_i') \beta(da_i') \biggr\rangle,
\end{eqnarray*}
where $\beta$ denotes the distribution of $a_1$.
Hence, we obtain
\begin{eqnarray*}
\hspace*{-5pt}&&\langle X_iY_i\rangle-\langle X_{i-1}Y_{i-1}\rangle\\
\hspace*{-5pt}&&\qquad =
\biggl\langle\int X_i(a_1,\ldots,a_{i-1},a_i')Y_i(a_1,\ldots,a_{i-1},a_i')
\beta (da_i') \biggr\rangle \\
\hspace*{-5pt}&&\qquad \quad{} -\biggl\langle\int X_{i}(a_1,\ldots,a_{i-1},a_i') \beta(da_i')\int
Y_{i}(a_1,\ldots,a_{i-1},a_i'') \beta(da_i'') \biggr\rangle
\\
\hspace*{-5pt}&&\qquad =
\biggl\langle\int\int\frac{1}{2}
\bigl(X_i(a_1,\ldots,a_{i-1},a_i')-X_i(a_1,\ldots,a_{i-1},a_i'') \bigr) \\
\hspace*{-5pt}&&\qquad \quad\hphantom{\biggl\langle\int\int} {} \times\bigl(Y_i(a_1,\ldots,a_{i-1},a_i')-Y_i(a_1,\ldots
,a_{i-1},a_i'') \bigr) \beta(da_i')
\beta(da_i'') \biggr\rangle\\
\hspace*{-5pt}&&\qquad \leq
\biggl\langle{\int\int\frac{1}{2}} \bigl(X_i(a_1,\ldots
,a_{i-1},a_i')-X_i(a_1,\ldots,a_{i-1},a_i'') \bigr)^2 \beta(da_i') \beta
(da_i'') \biggr\rangle^{1/2} \\
\hspace*{-5pt}&&\qquad \quad\!\! {} \times\biggl\langle{\int\int\frac{1}{2}} \bigl(Y_i(a_1,\ldots
,a_{i-1},a_i')-Y_i(a_1,\ldots,a_{i-1},a_i'') \bigr)^2 \beta(da_i') \beta
(da_i'') \biggr\rangle^{1/2}.
\end{eqnarray*}
We then conclude the proof of (\ref{eq:pr-covar}) as in the proof of
\cite{Gloria-Otto-09}, Lemma~2.3.

%s4.2 ###
\subsection{\texorpdfstring{Proof of Lemma \protect\ref{lem:ptwise-Green-decay}}{Proof of Lemma 4}}

We divide the proof in two main parts and deal with $|z|\leq\sqrt{T}$
and $|z|> \sqrt{T}$ separately.
The proof relies on the Harnack inequality on graphs.
We refer to Zhou \cite{Zhou-93} for $\Z^d$, and to Delmotte \cite
{Delmotte-97} for other graphs.
We recall here the easy part of Harnack's inequality (see \cite
{Delmotte-97}, Proposition~5.3,
or \cite{Zhou-93}, Proof of Theorem~3.3, (3.11)).
\begin{lemma}[(Harnack's inequality)]\label{lem:Harnack}
Let $a \in\calA$ and $R\gg1$. If $g\dvtx\Z^d\to\R^+$ satisfies
%
%e4.5 ###
\begin{equation}
-\nabla^*\cdot A\nabla g(x) \leq0
\end{equation}
in the annulus $\{R/2 < |x| \leq4R\}$ (i.e., $g$ is a nonnegative
subsolution), then
%
%e4.6 ###
\begin{equation}\label{eq:harnack}
\sup_{R<|x|\leq2R} g(x) \lesssim\biggl(R^{-d}\int_{R/2<|x|\leq4R}
g(x)^2\,dx \biggr)^{1/2}.
\end{equation}
\end{lemma}

\step{1} Proof of (\ref{eq:def_gT_d>2}) for $|x-y|\leq\sqrt{T}$.

Since $G_T$ satisfies
%
%e4.7 ###
\begin{equation}
-\nabla_x^*\cdot A\nabla_x G_T(x,y) = -T^{-1}G_T(x,y) \leq0
\end{equation}
for $|x-y|\gg1$, one may apply Lemma~\ref{lem:Harnack}.
For $R\gg1$, we then have
\[
\sup_{x:R<|x-y|\leq2R} G_T(x,y) \lesssim\biggl(R^{-d}\int
_{R/2<|x-y|\leq4R} G_T(x,y)^2\,dx \biggr)^{1/2}.
\]
Combined with \cite{Gloria-Otto-09}, Lemma~2.8, (2.21), for $q=2$
(which is uniform in $T>0$ and $y\in\Z^d$), this yields
\[
\sup_{R<|x-y|\leq2R} G_T(x,y) \lesssim R^{2-d},
\]
from which we deduce (\ref{eq:def_gT_d>2}) for $\sqrt{T} \geq|x-y|
\gg1$.
For $|x-y| \sim1$, we appeal to \cite{Gloria-Otto-09}, Proof of Lemma~2.8,
(4.4), with $R\sim1$ and $q=1$,
which yields $\sup_{|x-y|\leq R} G_T(x,y) \lesssim1$ by the
discrete $L^1-L^\infty$ estimate.

\step{2} Proof of (\ref{eq:def_gT_d=2}) for $|x-y|\leq\sqrt{T}$.

Let $N$ be a positive integer such that $2^N\sim\sqrt{T}$ and
$2^{-N}\sqrt{T}\gg1$.
For all $i\in\{1,\dots,N\}$, we first show that
\begin{eqnarray}\label{eq:ptwise-log-1-d=2}
&&\biggl(\bigl(2^{-i}\sqrt{T}\bigr)^{-2}\int_{2^{-i-1}\sqrt{T}<|x-y|\leq2^{-i+2}\sqrt
{T}} G_T(x,y)^2 \,dx \biggr)^{1/2}\nonumber\\[-8pt]\\[-8pt]
&&\qquad \lesssim i \sim\ln\biggl(\frac{\sqrt{T}}{1+2^{-i}\sqrt{T}}\biggr).\nonumber
\end{eqnarray}
Estimate \eqref{eq:ptwise-log-1-d=2} follows from the triangle
inequality and the BMO estimate of \cite{Gloria-Otto-09}, Lemma~2.8
(2.20), provided we show that
%
%e4.8 ###
\begin{eqnarray}\label{eq:ptwise-log-d=2}
\overline{G_T}_{\{|x-y|\leq2^{-i+2}\sqrt{T}\}} &\lesssim& i ,
\end{eqnarray}
where $\overline{G_T}_{\{|x-y|\leq2^{-i+2}\sqrt{T}\}} $ denotes the
average of $G_T(x,y)$ on the set $\{|x-y|\leq2^{-i+2}\sqrt{T}\}$.
By the triangle inequality and the BMO estimate of \cite
{Gloria-Otto-09}, Lemma~2.8, (2.20), we have
\begin{eqnarray*}
&&\overline{G_T}_{\{|x-y|\leq2^{-i+2}\sqrt{T}\}} \\[-2pt]
&&\qquad \leq \overline{G_T}_{\{|x-y|\leq2^{-i+3}\sqrt{T}\}} \\[-2pt]
& &\qquad \quad {} +2 \biggl(\frac{1}{|\{|x-y|\leq2^{-i+3}\sqrt{T}\}|}\\[-2pt]
&&\qquad \quad \hphantom{{} +2 \biggl(}
{}\times\int_{|x-y|\leq
2^{-i+3}\sqrt{T}}\bigl(G_T(x,y)-\overline{G_T}_{\{|x-y|\leq2^{-i+3}\sqrt
{T}\}}\bigr)^2\,dx \biggr)^{1/2}
\\[-2pt]
&&\qquad \leq \overline{G_T}_{\{|x-y|\leq2^{-i+3}\sqrt{T}\}} +C,
\end{eqnarray*}
where $C$ is a universal constant independent of $i$.
Combined with the estimate for $i=1$
\begin{eqnarray*}
{\overline{G_T}_{\{|x-y|\leq4\sqrt{T}\}} }&\lesssim& 1,
\end{eqnarray*}
which is a consequence of \cite{Gloria-Otto-09}, Lemma~2.8, (2.22),
this implies \eqref{eq:ptwise-log-d=2}
by induction.

We are now in position to prove (\ref{eq:def_gT_d=2}) for
$|x-y|\leq\sqrt{T}$.
Since $x\mapsto G_T(x,y)$ satisfies
\[
-\nabla_x^*\cdot A\nabla_xG_T(x,y) = -T^{-1}G_T(x,y) \leq0
\]
in the annulus $\{x, 2^{-i-1}\sqrt{T} < |x-y| \leq2^{-i+2}\sqrt{T}\}
$, Lemma~\ref{lem:Harnack} implies
\begin{eqnarray*}
&&\sup_{x:2^{-i}\sqrt{T}<|x-y|\leq2^{-i+1}\sqrt{T}} G_T(x,y)\\[-2pt]
&&\qquad \lesssim \biggl(\bigl(2^{-i}\sqrt{T}\bigr)^{-2}\int_{2^{-i-1}\sqrt{T}<|x-y|\leq
2^{-i+2}\sqrt{T}} G_T(x,y)^2 \,dx \biggr)^{1/2} \\[-2pt]
&&\qquad \lesssim \ln\biggl(\frac{\sqrt{T}}{1+2^{-i}\sqrt{T}}\biggr)
\end{eqnarray*}
using (\ref{eq:ptwise-log-1-d=2}) for $2^{-i}\sqrt{T}\gg1$.
For $|x-y|\leq R \sim1$, we appeal to (\ref{eq:ptwise-log-d=2}) and
to the discrete $L^1-L^\infty$ estimate
\[
G_T(x,y) \leq{R}^2 \overline{G_T}_{\{|x-y|\leq R\}} \lesssim
\ln T.
\]
This completes the proof of (\ref{eq:def_gT_d=2}) for $|x-y|\leq\sqrt{T}$.

\step{3} Proof of (\ref{eq:def_gT_d>2}) and (\ref{eq:def_gT_d=2}) for
$|x-y|>\sqrt{T}$.

Let $R\geq\sqrt{T}$.
Since $G_T$ satisfies
\[
-\nabla_x^* \cdot A \nabla_x G_T(x,y) = -T^{-1}G_T(x,y) \leq0
\qquad\mbox{for }|x-y|\geq1,
\]
Lemma~\ref{lem:Harnack} implies
\begin{eqnarray*}
\sup_{x:R<|x-y|\leq2R} G_T(x,y)
&\lesssim& \biggl(R^{-d}\int_{R/2<|x-y|\leq4R} G_T(x,y)^2\, dx
\biggr)^{1/2}.\vadjust{\goodbreak}
\end{eqnarray*}
Combined with \cite{Gloria-Otto-09}, Lemma~2.8, (2.23), for $q=2$ and
$r=k$, that is,
\begin{eqnarray*}
\int_{R/2<|x-y|\leq4R} G_T(x,y)^2 \,dx &\lesssim&R^{d+(2-d)2} \bigl(\sqrt
{T}R^{-1}\bigr)^{k},
\end{eqnarray*}
this yields the desired pointwise bound.

%s4.3 ###
\subsection{\texorpdfstring{Proof of Lemma \protect\ref{b-lem:dyad}}{Proof of Lemma 5}}

First note that by symmetry
\begin{eqnarray*}
\int_{|z|\leq|z-x|} h_T(z)h_T(z-x)\,dz&=&\int_{|z|\geq|z-x|}
h_T(z)h_T(z-x)\,dz \\[-2pt]
&\geq&\frac{1}{2}\int_{\Z^d} h_T(z)h_T(z-x)\,dz.
\end{eqnarray*}
Hence, it is enough to consider
\[
\int_{\Z^d} g_T(x)\int_{|z|\leq|z-x|}h_T(z)h_T(z-x)\,dz\,dx.
\]
In this proof, we essentially combine the pointwise decay of $g_T$ with
the results of \cite{Gloria-Otto-09}, Lemma~2.10,
that we recall here for the reader's convenience (see \cite
{Gloria-Otto-09}, Proof of
Lemma~2.10, Steps~1,~2 and 4):
There exists $\tilde R\sim1$ such that for all $R\geq\tilde R /2$,
\begin{eqnarray} \label{eq:Lemma-2.10-2}
&&\int_{R<|x|\leq2R}
\int_{|z|\leq|z-x|}h_T(z)h_T(z-x)\,dz\,dx\nonumber\\[-9pt]\\[-9pt]
&&\qquad \lesssim
\cases{
d=2,&\quad $R^2\max\bigl\{1,\ln\bigl(\sqrt{T}R^{-1}\bigr)\bigr\}$,\cr
d>2,&\quad $R^2$,
}\nonumber
\end{eqnarray}\vspace*{-12pt}
%e4.9 ###
\begin{eqnarray}
\label{eq:Lemma-2.10-1}
\int_{|x|\leq4\tilde R} \int_{|z|\leq|z-x|}h_T(z)h_T(z-x)\,dz\,dx
&\lesssim&\cases{
d=2,&\quad $\ln T$,\cr
d>2,&\quad 1.
}
\end{eqnarray}

In view of \eqref{eq:Lemma-2.10-1} and \eqref
{eq:Lemma-2.10-2}, it will be convenient to make\vspace*{1pt}
a dyadic decomposition of space. In order to also benefit from the
decay of $g_T(x)$ for $|x|\gg\sqrt{T}$, we
make the following decomposition of $\Z^d$:
%
%e4.12 ###
%e4.11 ###
%e4.10 ###
\begin{eqnarray}
\Z^d &=& \bigl\{|x| \leq2^{-I}\sqrt{T}\bigr\} \label{eq:decompZd-1}\\[-2pt]
&&{} \cup\bigsqcup_{i=-I,\dots,-1} \bigl\{2^i \sqrt{T} < |x| \leq
2^{i+1}\sqrt{T}\bigr\} \label{eq:decompZd-2}\\[-2pt]
&&{} \cup \bigsqcup_{i\in\N} \bigl\{2^i \sqrt{T} < |x| \leq
2^{i+1}\sqrt{T}\bigr\}, \label{eq:decompZd-3}
\end{eqnarray}
where $I$ is characterized by $2\tilde R < 2^{-I}\sqrt{T} \leq
4\tilde R$.

For the integral over the r.h.s. of \eqref{eq:decompZd-1},
we appeal to \eqref{eq:Lemma-2.10-1}
and to the definitions \eqref{b-eq:def_gT_d>2} and \eqref
{b-eq:def_gT_d=2} of $g_T(x)$ for $|x|\lesssim1$:
%
%e4.13 ###
\begin{eqnarray}\label{eq:lem-conv-1}
\qquad\quad\quad  \int_{ |x|\leq2^{-I} \sqrt{T} }g_T(x) \int_{|z| \leq
|z-x|}h_T(z)h_T(z-x)\,dz\,dx
&\lesssim& \cases{
d=2,& $\ln^2 T $,\cr
d>2,&1.
}\vadjust{\goodbreak}
\end{eqnarray}

For the integral over \eqref{eq:decompZd-3}, we use this
time \eqref{eq:Lemma-2.10-2} for $R\geq\sqrt{T}$
and the definitions \eqref{b-eq:def_gT_d>2} and \eqref
{b-eq:def_gT_d=2} of $g_T(x)$ for $|x|\geq\sqrt{T}$,
so that for all $i\in\N$ we have
\begin{eqnarray*}
&&\int_{ 2^i \sqrt{T} < |x| \leq2^{i+1} \sqrt{T} }g_T(x) \int_{|z|
\leq|z-x|}h_T(z)h_T(z-x)\,dz\,dx\\
&&\qquad \lesssim \bigl(2^i \sqrt{T}\bigr)^{2-d} (2^i)^{-3}\bigl(2^i\sqrt{T}\bigr)^2 \\
&&\qquad  =\sqrt{T}{}^{4-d}(2^i)^{1-d}.
\end{eqnarray*}
Summing this inequality on $i\in\N$ then yields the estimate
%
%e4.14 ###
\begin{eqnarray}\label{eq:lem-conv-2}
\int_{\sqrt{T} < |x|} g_T(x) \int_{|z| \leq|z-x|}h_T(z)h_T(z-x)\,dz\,dx
&\lesssim& \sqrt{T}^{4-d}.
\end{eqnarray}

We now deal with the integral over the last part \eqref
{eq:decompZd-2} of $\Z^d$.
To this aim, we combine \eqref{eq:Lemma-2.10-2} for $R\leq\sqrt{T}$
with the definitions \eqref{b-eq:def_gT_d>2}
and \eqref{b-eq:def_gT_d=2} of~$g_T(x)$ for $|x|\leq\sqrt{T}$.
In particular, for all $i\in\{-I,\dots,-1\}$, we have
\begin{eqnarray*}
&&\int_{ 2^i \sqrt{T} < |x| \leq2^{i+1} \sqrt{T} }g_T(x)
\int_{|z| \leq|z-x|}h_T(z)h_T(z-x)\,dz\,dx\\
&&\qquad \lesssim \cases{
d=2,&\quad $ \ln(2^{-i}) \bigl(2^i\sqrt{T}\bigr)^2\ln(2^{-i}) \sim i^2
\bigl(2^i\sqrt{T}\bigr)^2$, \cr
d>2,&\quad $\bigl(2^i\sqrt{T}\bigr)^2 \bigl(2^i\sqrt{T}\bigr)^{2-d} = \bigl(2^i\sqrt{T}\bigr)^{4-d}$.
}
\end{eqnarray*}
Summing this inequality over $i\in\{-I,\dots,-1\}$ and using that
$2^{I}\sim\sqrt{T}$ then yield
%
%e4.15 ###
\begin{eqnarray}\label{eq:lem-conv-3}
&&\int_{ 2^{-I} \sqrt{T} < |x| \leq\sqrt{T} }g_T(x) \int
_{|z| \leq|z-x|}h_T(z)h_T(z-x)\,dz\,dx\nonumber\\
&&\hspace*{10pt}\qquad \lesssim
\cases{
d=2,&\quad $\displaystyle  1+T\sum_{i=-I}^{-1} i^24^{i}$, \cr
d>2,&\quad $\displaystyle  1+\sqrt{T}^{4-d}\sum_{i=-I}^{-1}
(2^{4-d})^{i}$,
}
\\
&&\qquad  \stackrel{2^{I}\sim\sqrt{T}}{\lesssim}  \cases{
d=2,&\quad $ T$, \cr
d=3,&\quad $\sqrt{T}$, \cr
d=4,&\quad $\ln T$, \cr
d>4,&\quad 1.
}\nonumber
\end{eqnarray}

The combination of \eqref{eq:lem-conv-1}, \eqref{eq:lem-conv-2} and \eqref{eq:lem-conv-3} finally proves \eqref{b-eq:dyad}.

%%%%%%%%%%%%%%%%%%%%%%%%%%%%%%%%%%%%%%%%%%%%%%%%%%%%%%%%
%

% imsref loaded by dianan, 2011-02-17 15:06:23
%

\printaddresses

\end{document}